\magnification=1200

\input eplain

\input epsf

\newcount\secno
\newcount\thmno
\def\here{\the\secno.\the\thmno}

\newcount\figno

\def\sec#1{\par\vskip1cm
\advance \secno by 1
\centerline{{\bf \the\secno\ #1}}
\thmno=0
\par\bigskip\noindent}

\long\def\thm#1#2{\par\bigskip\noindent
\advance\thmno by 1
{\bf #1 \here.} {\it #2}
\smallskip}

\long\def\thmnamed#1#2#3{\par\bigskip\noindent
\advance\thmno by 1
{\bf #1 \here\ (#2).} {\it #3}
\smallskip}

\def\rmk#1{\par\medskip\noindent
\advance\thmno by 1
{\bf #1 \here. }}

\def\proof{\noindent{\it Proof: }}

\def\figure#1{\global\advance\figno by 1
\bigskip
\centerline{Figure~\the\figno: #1}}

\def\eps{\varepsilon}
\def\pr{{\rm pr}}

\secno=-1

\everyfootnote={\leftskip=1cm \parindent=0cm 
                \rm}


\def\id{{\rm Id}}
\def\r{{\bf R}}
\def\z{{\bf Z}}
\def\im{{\rm im}}
\def\cstr{contact structure}
\def\ot{overtwisted}
\def\ocstr{overtwisted \cstr}
\def\difmm{diffeomorphism}
\def\cmm{contactomorphism}    
\def\lk{Legendrian knot}
\def\leg{Legendrian}

\def\cont{{\it Cont}}
\def\heq{homotopy equivalen}
\def\nbd{neighborhood}
\def\qed{{\hfill \boxit{}}}

\def\fott{front on the torus}
\def\fotts{fronts on the torus}
\def\cusp#1{cusp_{#1}}
\def\cross#1{cross_{#1}}
\def\tbi{Thurston-Bennequin invariant}
\def\consum{connected sum}

\def\g#1#2{\Gamma_{#1,#2}}


\def\be{\cite{bennequin}}
\def\dy{\cite{dymara:masters}}
\def\futa{\cite{fuchs-tabachnikov}}

\font\ec=ecrm10
\def\c{\hskip.27em{\ec \char'14}\hskip-.54em}

\font\scaps=cmcsc10

\centerline{\bf LEGENDRIAN KNOTS IN OVERTWISTED CONTACT STRUCTURES}

\bigskip
\centerline{Katarzyna Dymara\footnote{*}{Partially supported by KBN 
grant 2 P03A 017 25.}}
\centerline{Institute of Mathematics}
\centerline{Wroc\l aw University}
\bigskip
\centerline{October 2004}

\vskip1cm

{
\narrower

\centerline{\scaps Abstract}
{
We prove that two Legendrian knots in a contact
structure which is trivializable as a plane bundle
are Legendrian isotopic provided that
(1) they are isotopic as framed knots, (2) they have the same
rotation number with respect to some parallelization of the contact 
structure, and (3) there is an overtwisted disk disjoint with both 
knots. (For zero-homologous knots the condition (1) reads as: (1a)
they are isotopic as topological knots, and (1b) they have
the same Thurston-Bennequin invariant.) 
Then we discuss the situation when condition (3) is not fulfilled,
in particular that of non-loose Legendrian knots.}

}

\sec{Introduction}

The present paper deals with \lk s in 
\ocstr s.
In \dy, the author proved (using Eliashberg's and Gray's 
theorems recalled here as 
\refs{thm-eliashberg} 
and \refn{thm-grayrel} respectively)
that two \lk s in an \ocstr\ on $S^3$
disjoint with a fixed \ot\ disk
are \leg\ isotopic if and only if
they represent the same topological knot type and have equal values
of invariants $rot$ and $tb$ (see \ref{def-invts}).
The proof relies on homotopy considerations 
specific for $S^3$, which makes it difficult to generalize 
to \ocstr s on other manifolds.

In this paper we want to give a proof by direct manipulation
on fronts, with the help of the idea (introduced by Fuchs and 
Tabachnikov in \futa) of following an arbitrary isotopy of knots by a 
\leg\ isotopy adding zig-zags whenever necessary. 
Because of those additional zig-zags, this method does not 
yield actual \leg\ isotopy, but only equivalence in the
Grothendieck group of \leg\ knots (with the connected summation as 
the group operation).
As it turns out, in the presence of an \ot\ disk 
an operation functionally identical to adding zig-zags can 
be performed by \leg\ isotopy.
There is a clear advantage to this approach: 
one can fairly easily construct actual isotopies of \lk s
following the steps of the proof of \ref{thm-main}
(the isotopies given as examples in \ref{sec-lks} have been
found this way), which would be rather problematic with
the methods of \cite{dymara:masters}.
On the other hand, comparing the two proofs suggests
some sort of connection between Fuchs and Tabachnikov's idea of
\lk\ stabilization and the notion of \ot ness, which may be
interesting in itself.

The paper is organized in the following way.
\ref{sec-cstr} recalls
definitions and known facts on \cstr s, and fixes
notation to be used throughout the paper.
In \ref{sec-lks} we turn attention to \lk, again
restating known results and
adapting them to the situation of an
\ocstr\ on $S^3$; this process continues throughout 
\ref{sec-invts}, but with respect to the classical invariants.
\ref{sec-mainthm} formulates and proves the main theorem
(\ref{thm-main}).
Finally, in (somewhat eclectic) \ref{sec-zoo} we display 
several supposedly interesting examples of \lk s,
including many non-loose ones,
and give a few classification results.

\sec{Contact Structures}
\definexref{sec-cstr}{\the\secno}{section}

\noindent
A {\it \cstr} on a 3--dimensional manifold is a field of planes 
defined (at least locally) as the kernel of a 1--form 
$\alpha$ such that $\alpha\wedge d\alpha$ nowhere vanishes.


\rmk{Example} \definexref{ex-std}{\here}{exple} 
The \cstr\ $\zeta=\ker(\alpha)$, where $\alpha=dz-x\,dy$, 
is called the standard \cstr\ on $\r^3$.

\rmk{Example} \definexref{ex-stdprim}{\here}{exple}
The \cstr\ $\zeta'=\ker(\alpha')$, 
where $\alpha'=\cos x\,dz-\sin x\,dy$,
is isomorphic to $\zeta$ of \ref{ex-std}
via the \difmm\ 
$$f_1:\r^3\to\r^3:(x,y,z)\mapsto(x\, y\,\cos x+z\,\sin x,
\ x\,y\,\cos x-y\,\sin x+x\,z\,\sin x + z\,\cos x).$$


\midinsert
\centerline{\epsffile{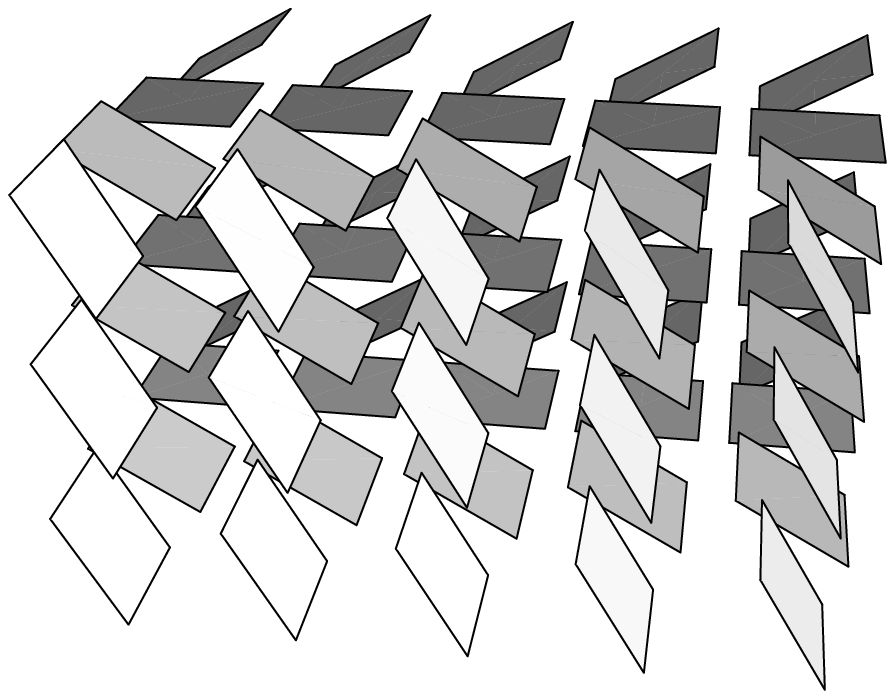}}
\figure{The standard \cstr\ $\zeta$.}
\endinsert

\bigskip

We say that two \cstr s $\zeta=\ker\alpha$ and $\xi=\ker\beta$ on a 
manifold $M$ are {\it isomorphic} if 
there exists a \difmm\ $f:M\to M$ such that $f^*(\alpha)=t\beta$ 
for some non-zero function $t:M\to\r$.

\smallskip
In general, it is true that all \cstr s look alike {\it locally}
(but not globally, cf.~\ref{ex-ot});
to be precise, in contact topology holds 
an analogue of the well-known (in symplectic setting) Darboux's Theorem.


\thmnamed{Theorem}{Darboux's Theorem}{In any contact manifold 
any point has a \nbd\ $U$ 
contactomorphic to the standard \cstr\ $\zeta$ 
{\rm (see \ref{ex-std})}
on $\r^3$.}


(For both the symplectic Darboux's Theorem and its contact version, 
consult for 
example Arnold's book \cite{arnold}.)

\rmk{Example} \definexref{ex-ot}{\here}{exple}
The \cstr\ $\xi$ on $\r^3$, defined as the kernel of the 1--form
(in cylindrical coordinates) $\beta=\cos r\,dz-r\,\sin r\,d\theta$, 
is not isomorphic to the \cstr\ $\zeta$
of \ref{ex-std}. (Essentially, what 
distinguishes $\xi$ from $\zeta$ is 
its property of being {\it \ot}, to be explained in \ref{def-ot}
and the discussion following it.)

\smallskip
Apart from the relation of isomorphism, 
there is another (equally natural, one would say, though possibly weaker) 
equivalence relation on the space of \cstr s on a given manifold, namely 
that of isotopy (homotopy through \cstr s). 

\rmk{Definition}
We say that two \cstr s $\zeta_0$ and $\zeta_1$ are 
{\it isotopic} if there exist a smooth family 
of \cstr s $\{\zeta_t, t\in[0,1]\}$.

\medskip
An easy argument shows that on $\r^3$ all \cstr s are isotopic. 
For $\xi$ a \cstr, let $\eps$ be such that $\xi|_B$, where $B=B_\eps(0)$, 
is isomorphic to the standard
\cstr\ $\zeta$ (see \ref{ex-std}). 
Then the family of pull-backs of $\xi$ by 
maps $f_t$, which gradually shrink the entire $\r^3$ 
so that finally it falls into the ball $B$,
provides an isotopy between $\xi$ and $\zeta$: for 
\vskip-3mm
$$f_t:\r^3\to\r^3:x\mapsto\left({\eps \over t\|x\|+1}\right)x 
\qquad\hbox{and}\qquad
\xi_t=f_t^*(\xi)$$
we have $f_0=\id$, so that $\xi_0=\xi$, and  $f_1(\r^3)=B$. Therefore 
$\xi_1$ is isomorphic to the standard $\zeta$.

The situation is quite different in the case of closed manifold: 
Gray proved in \cite{gray} that the relation of isotopy is 
as strong there as that of isomorphism.


\thmnamed{Theorem}{Gray's Theorem\numberedfootnote{Actually, 
the original Gray's formulation goes along more telling 
and easy to remember lines:
{\it Families of contact structures are locally trivial} (Theorem 5.2.1 
in \cite{gray}).}}{Given a smooth family $\{\zeta_t, t\in D^n\}$ of
contact structures on a closed manifold $M$ and a point $t_0\in D^n$,
there exists a family $\{\phi_t, t\in D^n\}$ of \difmm s $\phi_t:M\to M$ 
such that
for all $t$, $\phi^*\zeta_t=\zeta_{t_0}$.}
\definexref{thm-gray1}{\here}{thm}


\medskip
Throughout the paper, we will happily make use of 
this theorem, 
calling two \cstr s ``isomorphic'' (or simply ``the same'') as soon as an 
isotopy 
between them has been established. 

Gray's proof (which uses a vector field whose flow consists of 
the desired  \difmm s, constructed locally 
and glued together by means of a partition of unity) 
can be applied---without essential changes---to a relative situation, 
either in the sense of considering \cstr s on a manifold modulo
a compact set, or fixing the \cstr\ along a (contractible)
subset of the parameter space. Here we will simply reformulate
Gray's theorem in a way as general as it is needed.

\thmnamed{Theorem}{Gray's Theorem, relative version}{Let 
$M$ and $\zeta_t$ be as in \ref{thm-gray1}.
Assume that $\zeta_t|_A=\zeta_{t_0}|_A$ for a compact set $A\subset M$ 
and for all $t\in D^n$. Moreover, let $\zeta_t=\zeta_{t_0}$ for  all  
$t\in D'$, where $D'$ is a contractible subset of $D^n$. 
Then there exists a family $\{\phi_t, t\in D^n\}$ of \difmm s 
$\phi_t:M\to M$ such that for all $t$, $\phi^*\zeta_t=\zeta_{t_0}$, 
for all $t\in D^n$ $\phi_t|_A=\id_A$ and for all $t\in D'$ $\phi_t=\id_M$.}
\definexref{thm-grayrel}{\here}{thm}

This is just one of the reasons why in the present paper 
we prefer to discuss \cstr s on a
closed manifold. Since all the basic examples described above happened 
to be \cstr s on $\r^3$ (due to 
the ease with which to set a coordinate system on the euclidean space), 
we will make up for this one-sided approach by providing a spectrum 
of examples sitting on the 3-sphere.

\bigskip
In order to accompany these examples by an illustrative figure, 
let us employ the following way of looking at $S^3$.

\midinsert
\centerline{\epsffile{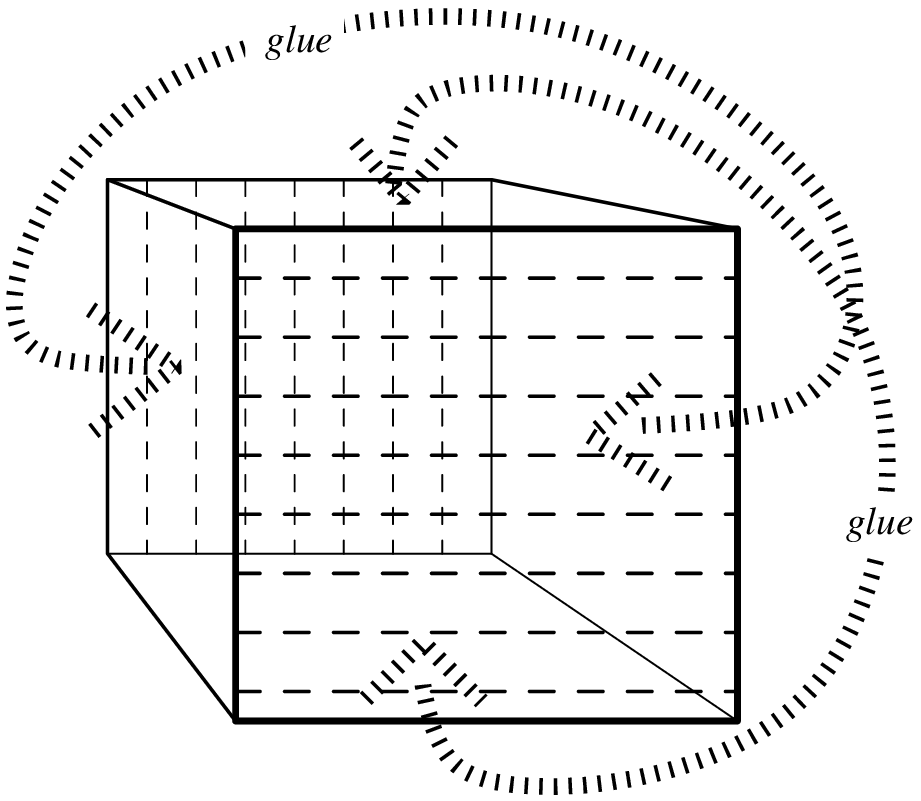}}
\figure{Gluing of the cube (\ref{cube}).}
\definexref{fig-cube}{\the\figno}{figure}
\bigskip
\endinsert

\rmk{Construction} \definexref{cube}{\here}{construction}
Consider the cube  ${\cal C}=[0,1]\times[0,1]\times[0,1]$. 
Identify two points $(x,y,z)$ and $(x',y',z')$ in $\cal C$
if any of the following condition holds:\numberedfootnote{In 
more informal language we would say: glue together 
the left and right faces of the cube, 
and also the top and bottom faces; then contract to points all 
horizontal lines on the front face, and also all vertical lines 
on the back face. Cf.~\ref{fig-cube}.}
\item{(a)} $x=x', y=y', z=0, z'=1$, or
\item{(b)} $x=x', y=0, y'=1, z=z'$, or
\item{(c)} $x=x'=0, z=z'$, or
\item{(d)} $x=x'=1, y=y'$.

Then the quotient is homeomorphic to $S^3$.
We will denote the projection map ${\cal C}\to S^3$
by $\pr_{\cal C}$.

\medskip

This construction provides a system of coordinates on the 3-sphere, 
which can be used to write explicitly 1-forms defining certain \cstr s.

\rmk{Example} \definexref{ex-zetan}{\here}{exple}
For any integer $n\geq 0$, let the \cstr\ $\zeta_n$ on $S^3$ 
be the kernel of the 1-form 
$\alpha_n=\cos\bigl( (n+\frac1/2)\pi x\bigr)\,dz
-\sin\bigl( (n+\frac1/2) \pi x\bigr)\,dy$.
\smallskip


Note that for each $n$ the form $\alpha_n=dz$ on the front face 
$\{x=0\}$ and $\alpha_n=\pm dy$ on the back face $\{x=1\}$.
Thus all $\zeta_n$'s are indeed well-defined \cstr s on $S^3$.
Now let us have a closer look at one of them.

\rmk{Example}
\definexref{ex-stdonsphere}{\here}{exple}
The \cstr\ $\zeta_0$ (defined as above with $n=0$)
is known as the standard \cstr\ on $S^3$. 
It can be described independently of \ref{cube}
as the field of planes perpendicular 
to the fibers of the Hopf bundle. 
Removing a point from $(S^3,\zeta_0)$
yields a \cstr\ on $\r^3$ isomorphic to $\zeta$ from 
\ref{ex-std}. 

\medskip
Each of the \cstr s described in 
\ref{ex-zetan} belongs to one of exactly three different 
isomorphism classes. It is easy to notice that, 
considered as plane fields, all $\zeta_n$'s
with $n$ odd realize one homotopy type, those with $n$ even---another one.
But in addition to that, $\zeta_0$ can be proved to be different
from $\zeta_2$ ($\simeq\zeta_4\simeq\ldots$) using the notion of 
\ot ness (introduced in \cite{bennequin}).


\rmk {Definition} \definexref{def-ot}{\here}{defin}
Let $R=[0,\frac2/3]\times[0,1]\times\{z_0\}\subset{\cal C}$
for some $z_0$,
$\Delta=\pr_{\cal C}(R)$ in $S^3$ with the \cstr\ $\zeta_1$  of \ref{ex-zetan}.
A 2-dimensional disk $D$ in a contact manifold $(M,\xi)$ is an 
{\it \ot} disk if it has a \nbd\ contactomorphic to a \nbd\ of $\Delta$.
A \cstr\ (or contact manifold) is called {\it \ot} if there is 
an \ot\ disk in it.
\medskip


\rmk{Example} \definexref{ex-otdisk}{\here}{exple}
Let $R'=[\frac1/3,1]\times\{y_0\}\times[0,1]\subset {\cal C}$
for some $y_0$. 
The disk $\pr_{\cal C}(R')$ is another \ot\ disk in $(S^3,\zeta_1)$.

\smallskip
\rmk{Example}
For each $z_0$, the disk $\Delta_{z_0}=\{(r,\theta,z) \mid
r\leq\pi, z=z_0\}$ is an \ot\ disk in $(\r^3,\xi_0)$ for $\xi_0$ 
defined as in \ref{ex-ot}.

\smallskip
\rmk{Definition}
A \cstr\ which is not \ot\ is called {\it tight}.

\smallskip
\rmk{Example}
The \cstr\ $\zeta$ (and $\zeta'$) on $\r^3$ described in
\refs{ex-std} and \refn{ex-stdprim}, as well as $\zeta_0$ on $S^3$ 
from \ref{ex-stdonsphere}, 
are tight.

\medskip
The tight-vs-overtwisted dichotomy has unexpectedly deep 
consequences. Results obtained and methods useful in
studying \ot\ contact manifolds 
have very little in common with those appropriate
for tight ones. The classification results 
for tight \cstr s are rather mysterious: 
not universal but proved on a manifold-by-manifold basis 
(\cite{honda:classification1}, \cite{honda:classification2}, 
\cite{etnyre:tight_lens}),
often surprising (\cite{non-existence}) and refer to a
rich hierarchy of notions of fillability (cf.~the discussion of
Question~1 in \cite{problems}).
In the case of \ocstr s virtually all that has ever been said 
about classification is contained in one powerful 
(though deceivingly simple in formulation) 
theorem of Eliashberg \cite{eliashberg:classification}.


\thmnamed{Theorem}{Eliashberg's Theorem}{For $\Delta$ a 2-dimensional disk 
in an arbitrary 3-manifold $M$,
let $\xi_\Delta$ be a \cstr\ on a \nbd\ of $\Delta$ 
for which $\Delta$ is an \ot\ disk.
Denote by $\cont_\Delta(M)$ the space of \cstr s on $M$
which coincide with $\xi_\Delta$ on a \nbd\ of $\Delta$.
Moreover, let ${\it Distr}_\Delta(M)$ be the space of all plane
distributions on $M$ which coincide with $\xi_\Delta$ on a \nbd\ of $\Delta$.
Then the natural embedding
$\cont_\Delta(M)\to{\it Distr}_\Delta(M)$ is a weak \heq ce. }
\definexref{thm-eliashberg}{\here}{thm}


\smallskip
On $\pi_0$ level this means that isotopy classes of 
\cstr s \ot\ along a fixed disk (which, since embeddings of a disk
into a connected manifold are all isotopic, actually exhaust 
all isotopy classes of \ocstr s) remain in a one-to-one correspondence
with homotopy classes of plane fields.

For $S^3$, which is our main object of interest, 
co-oriented\numberedfootnote{Of course, all plane fields on $S^3$
(or any manifold with trivial $H_1$) are co-orientable; we consider
them as pre-equipped with one of the two possible co-orientations
simply as a matter of convenience.}
plane fields can be identified with maps into $S^2$ by means of 
fixing a parallelization (a trivialization of the tangent bundle).
Therefore, the space of \ocstr s is parameterized by the set
$\pi_0\bigl({\it Map}(S^3{\to}S^2)\bigr)=\pi_3S^2\simeq{\bf Z}$.

\sec{Legendrian Knots}
\definexref{sec-lks}{\the\secno}{section}

Now we turn attention to \lk s.
In this section (with its multitude of figures)
we introduce a device, analogous to front diagrams in the standard
tight \cstr, but suitable for studying \lk s in in \ocstr s on $S^3$
(\ref{fronts-on-torus}), display several examples,
and prove a(n appropriate version of) theorem about Reidemeister moves.

A curve in a contact 3-manifold tangent to the contact planes at each point
is called {\it \leg}. The term {\it \lk} is usually applied
to a \leg\ curve which is a knot, i.e.~an embedding of the 
circle.\numberedfootnote{In \futa, 
Fuchs and Tabachnikov refer to such a curve 
(in the standard tight \cstr\ on $\r^3$) as a {\it short} \lk, 
with the term {\it long \lk s} reserved for \leg\ embeddings of $\r$
coinciding with the $y$ axis outside of a compact set;
one can broaden the notion of a long \lk\ 
to include any proper \leg\ embedding of $\r$ 
(with fixed behavior in the infinity) into an open
contact manifold. All the results of this paper are 
valid for long as well as short \lk s, unless otherwise specified.}

We say that two \lk s are {\it \leg\ isotopic} if
they are isotopic via \lk s; this is an equivalence relation.
The words ``\lk'' are often used to mean 
a \lk\ {\it type}, i.e.~a class of \leg\ isotopy.

We will also make extensive use of {\it knot diagrams},
i.e.~projections of knots on a plane (or a surface) in general position.

Projections of \lk s in $(\r^3,\zeta)$
on the $yz$ plane yield their so called {\it front diagram}.
Those are plane curves which are never vertical 
and may have---despite the knot being in general position---cusp 
singularities. Fronts are easier to manipulate, since there is no
metric condition involved, and provide a comfortable tool 
for studying \leg\ isotopies.
Two \lk s are \leg\ isotopic if and only if there is a homotopy 
between their front diagrams which is an isotopy except for 
finitely many points, where one of the {\it \leg\ Reidemeister moves}
(recalled on \ref{fig-moveslocal}) occurs \cite{jacek}.

In the present paper we will employ a device similar to front diagrams
in $\r^3$, but better suited to the contact manifold of our interest.

\rmk{Construction} \definexref{fronts-on-torus}{\here}{construction}
({\it Fronts on the torus, fronts on the cylinder.})
Consider $S^3$ with the \cstr\ $\zeta_1$ as in \ref{ex-zetan}. 
Recall \ref{cube}; the images of the front and back faces
are two linked circles. 
A \lk\ $K$ in general position is disjoint with those circles, 
i.e.~$\pr_{\cal C}^{-1}(K)$ lies in ${\cal C}_0=(0,1)\times[0,1]\times[0,1]$.
Let $p:{\cal C}_0\to\{0\}\times[0,1]\times[0,1]$ 
be the orthogonal projection; since $\pr_{\cal C}$ is injective on ${\cal C}_0$, 
$P=\pr_{\cal C}\circ p\circ\pr_{\cal C}^{-1}$ is well-defined and 
$P(K)$ is a curve on the torus $T^2=\pr_{\cal C}(\{0\}\times[0,1]\times[0,1])$,
smooth except for possible cusp singularities.
We will refer to $P(K)$ as the {\it front diagram} of $K$ {\it 
on the torus}.

Note that $\r^3$ with the \cstr\ $\xi$ of \ref{ex-ot}
is contactomorphic to the open subset 
$\pr_{\cal C}\bigl([0,1)\times[0,1]\times(0,1)\bigr)$.
Thus a \lk\ in $(\r^3,\xi)$ can be diagrammatically represented by
a \fott\ never crossing a fixed meridian of the torus
(i.e.~the image of the top and bottom faces of the cube);
such a curve may be called a {\it front on the cylinder}.

Moreover, let $\Delta$ be an \ot\ disk in a contact manifold $(M,\zeta)$.
A small \nbd\ of $\Delta$ is contactomorphic to a \nbd\ of the standard
\ot\ disk in $(S^3,\zeta_1)$ 
(by the standard \ot\ disk in $(S^3,\zeta_1)$ we understand 
$pr_{\cal C}^{-1}(R)$ for 
$R=[0,\frac2/3]\times [0,1]\times\{z_0\}\subset 
{\cal C}$; cf.~\ref{ex-otdisk}).
Therefore \leg\ curves in this \nbd\ can also be represented 
by fronts on the cylinder. 

\midinsert
\centerline{\epsffile{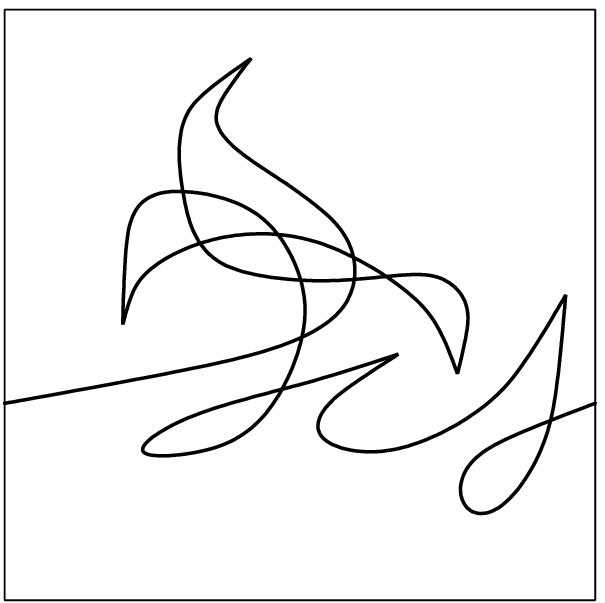}}
\figure{A front on the torus (which can be 
interpreted as a front on the cylinder).}
\definexref{fig-afront}{\the\figno}{figure}
\endinsert

\medskip
It is not the case that every curve on the torus is the projection of
a \lk. Knot diagrams on the torus satisfy a condition analogous to---though 
more complicated than---that of front diagrams being never vertical.

For $n\geq 1$, a curve of any slope can be the projection of a \leg\ curve
in $\zeta_n$.
Moreover, for some slopes there is an ambiguity 
in determining the $x$ coordinate of the original curve. 
Indeed, the \cstr\ $\zeta_1$ looks exactly the same 
at $(x,y,z)$ as at $(x+\frac{2}/{3},y,z)$;
and generally: $\zeta_n$ looks exactly the same 
at $(x,y,z)$ as at $(x+2/(2n+1),y,z)$.
Thus locally any curve in the torus lifts to
$n$ or $n+1$ different curves in $(S^3,\zeta_n)$. 
However, once a starting point has been chosen,
the \leg\ lift has to be continued in a unique way;
and it may happen that a particular curve cannot be lifted. 
Below we formulate the necessary and sufficient conditions
for a curve on the torus to be the projection of a \lk.

\thm{Proposition}{Let $k:S^1\to T^2$ be a curve in the torus, 
$g_k:S^1\to RP^1$ its tangent Gauss map.\numberedfootnote{We 
take it to have the codomain in $RP^1$ rather than $S^1$ so that
we don't have to fix an orientation of $k$.} 
Let $p:\r\to RP^1$ (of period $\pi$) be the universal covering map
of $RP^1$.
The curve $k$ is the projection of a \lk\ in $\zeta_n$ if and only if
\item{$(i)$} $g_k$ lifts to a map $\tilde g_k:S^1\to\r$, 
where $g_k=p\circ \tilde g_k$; and 
\item{$(ii)$} $\im(\tilde g_k)\subset \big(\ell\pi,(\ell+n+\frac1/2)\pi\big)$
for some $\ell\in {\bf Z}$.}
\definexref{propo-lift}{\here}{propo}

\proof
Assume that the curve $k$ is the projection of a \lk\ $K$, 
where $K(t)=\bigl(x_K(t),y_K(t),z_K(t)\bigr)$. 
The Gauss map $g_k$ must agree with the $x$ coordinate of $K$,
i.e.~as $\tilde g_k$ we can choose $\frac3/2\pi x_K(t)$.
But $\im(x_K)\subset (0,1)$, so the condition $(ii)$
is satisfied.

On the other hand, when
$\im(\tilde g_k)\subset 
\bigl(\ell\pi, (\ell+n+\frac1/2)\pi\bigr)$, 
take $x_K(t)=\tilde g_K(t) - \ell\pi$.
\qed

\medskip
Even for $n=1$ we can notice that the intervals 
$\bigl(\ell\pi, (\ell+\frac3/2)\pi\bigr)$
are not pairwise disjoint.
This is connected with the ambiguity in determining the
$x$ coordinate discussed above.
If the image of $\tilde g_k$ is contained in the interval
$(0,\frac{$\pi$}/2)$---as it is, for instance, for the flying saucer
of \ref{fig-ufo}---then $k$ is the diagram of each of two distinct
\lk s in $\zeta_1$, one of which lies entirely in the front $\frac1/3$ of the 
cube, and the other one in the back $\frac1/3$.
However, the difference of the values of $\tilde g_k$ (or the 
$x$ coordinates) at two given points on a knot is always
uniquely determined. This observation will play an important
role in \ref{propo-reidemoves}. 

\rmk{Example} \definexref{ex-gammas}{\here}{exple}
A curve on the torus of constant slope $m/n$ (see \ref{fig-gammas})
is the diagram of a \lk\ in $\zeta_1$, which we will denote $\g mn$.
If $mn\leq 0$, then $\g mn$ is
uniquely determined; otherwise, there are two possibilities.
For example,
$\Gamma_{0,1}$ and $\Gamma_{1,0}$ have diagrams consisting of
(respectively) a horizontal and a vertical line segment.
They are boundaries of \ot\ disks described in \ref{def-ot}
and \ref{ex-otdisk}. 

\midinsert
\centerline{\epsffile{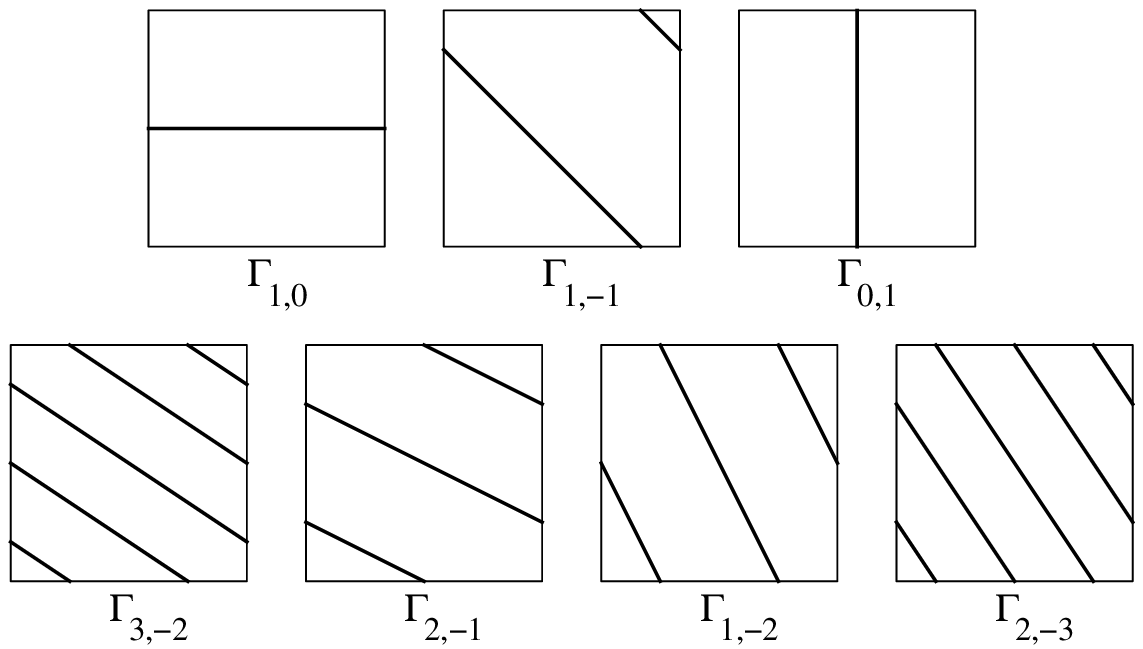}}
\smallskip
\centerline{(a) the uniquely determined $\Gamma_{m,n}$'s}
\vskip14mm
\centerline{\epsffile{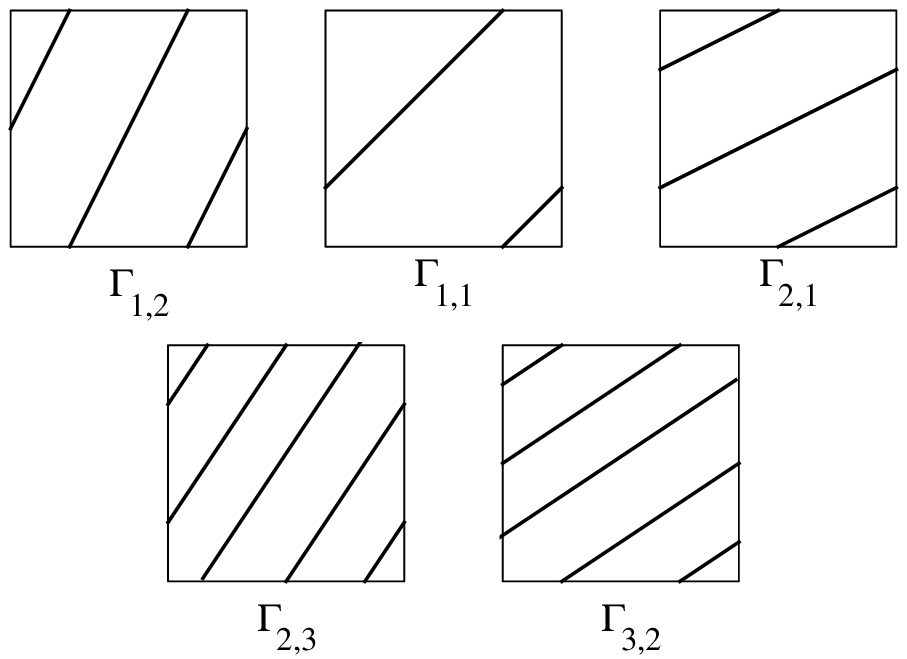}}
\smallskip
\centerline{(b) the ambiguous $\Gamma_{m,n}$'s}
\bigskip
\figure{The knots $\Gamma_{m,n}$.}
\definexref{fig-gammas}{\the\figno}{figure}
\endinsert

\medskip
Fronts on the torus may be used in studying \leg\
isotopies of knots, in a similar way as \hbox{\'Swi\c at}\-kow\-ski
in \cite{jacek} used the usual front diagrams for knots
in the standard tight \cstr\ on $\r^3$. 
The following result is analogous to his Theorem~B.

\thm{Proposition}{Two fronts on the torus
are projections of \leg\ isotopic knots if and only if
we can pass from one to the other by a finite sequence 
of moves of the following types:
\item{\rm (0)} an isotopy of the fronts;
\item{\rm (I)} the type I Reidemeister move (the dovetail move), see \ref{fig-moves1};
\item{\rm (II)} the type II Reidemeister move in two versions:
\itemitem{\rm (a)} in the presence of a cusp, see \ref{fig-moves2a}, and
\itemitem{\rm (b)} without a cusp (see \ref{fig-moves2b}), 
but under the assumption that the $x$ coordinates of both strands
differ by at least $\pi$;\numberedfootnote{This condition has been
symbolically depicted in \ref{fig-moves2b} by joining the strands with a dotted line
which makes a U-turn before returning to the point where the move is performed;
of course, the actual diagram may look different.}
\item{\rm (III)} the type III Reidemeister move, see \ref{fig-moves3}; and
\item{\rm (IV)} the move changing the homotopy class of the diagram 
(see \ref{fig-moves4}) by passing through
\itemitem{\rm (a)} the front face of the cube, or
\itemitem{\rm (b)} the back face.}
\definexref{propo-reidemoves}{\here}{propo}

\midinsert
\centerline{\epsffile{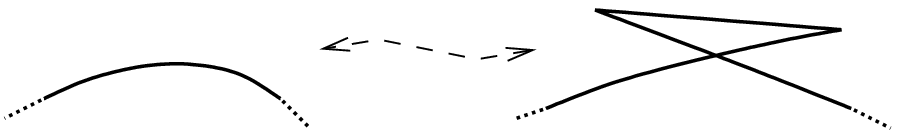}}
\centerline{The type I Reidemeister move}
\centerline{\epsffile{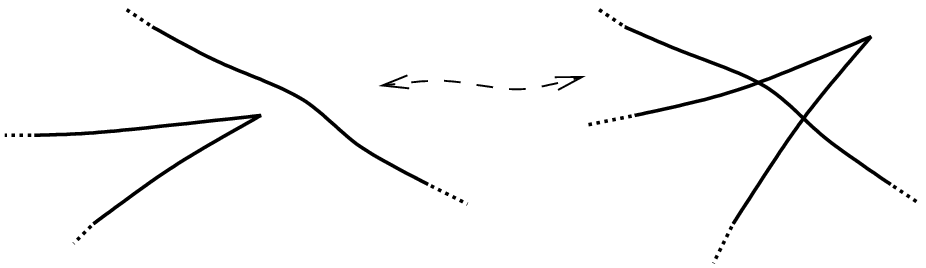}}
\centerline{The type II Reidemeister move (involving a cusp)}
\centerline{\epsffile{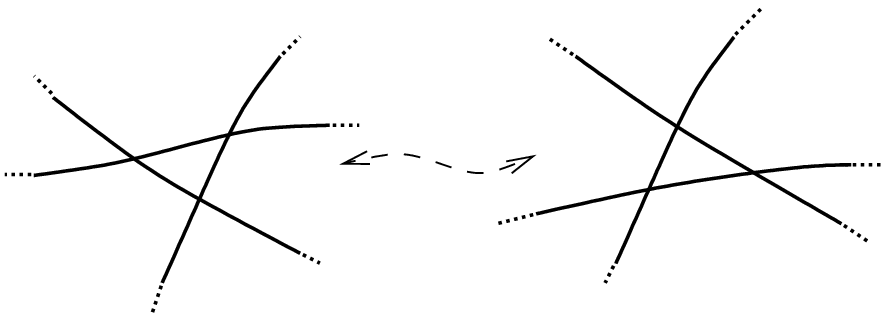}}
\centerline{The type III Reidemeister move}
\figure{The local Legendrian Reidemeister moves (cf.~\cite{jacek}).}
\definexref{fig-moves1}{\the\figno(I)}{fig}
\definexref{fig-moves2a}{\the\figno(II)}{fig}
\definexref{fig-moves3}{\the\figno(III)}{fig}
\definexref{fig-moveslocal}{\the\figno}{figure}
\endinsert

\midinsert
\centerline{\epsffile{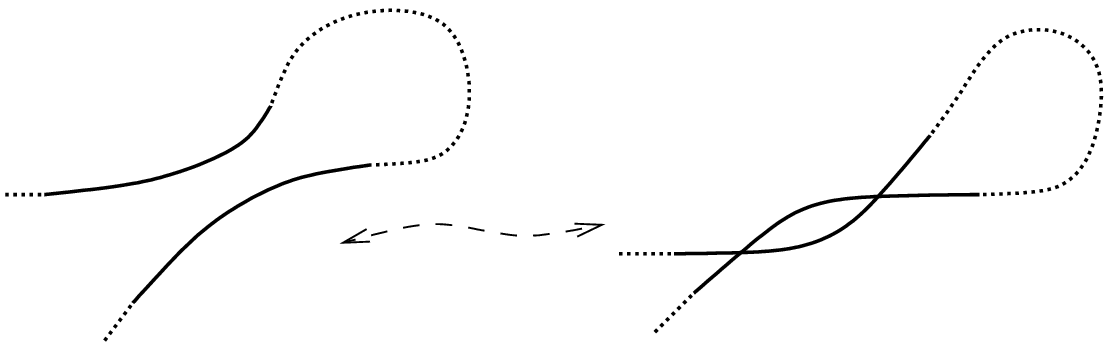}}
\figure{The type II Reidemeister move without a cusp.}
\definexref{fig-moves2b}{\the\figno}{fig}
\endinsert

\midinsert
\null\bigskip
\centerline{\epsffile{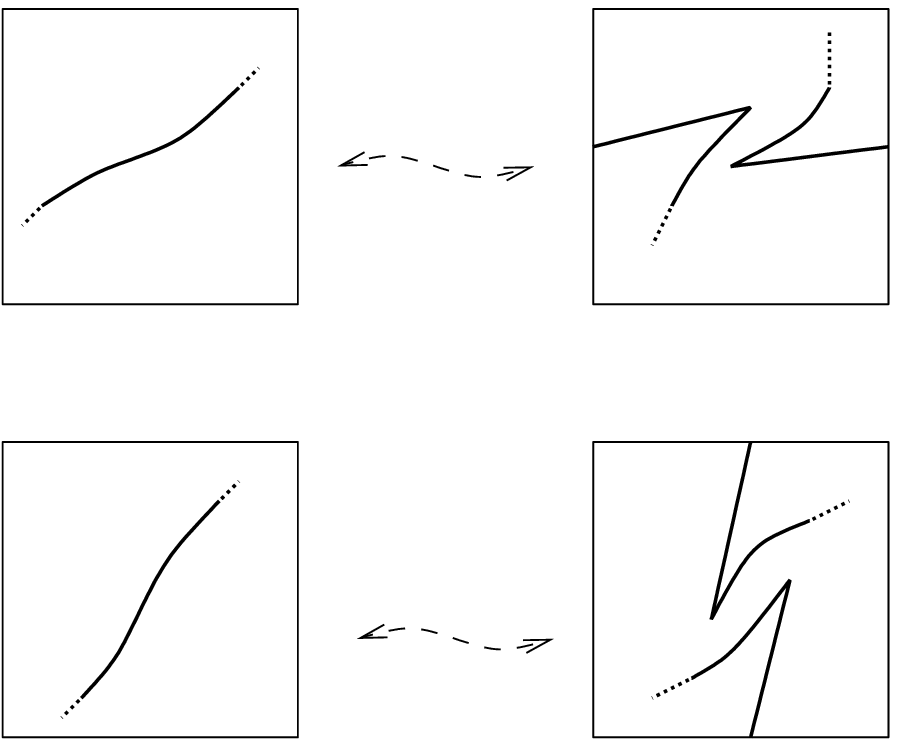}}
\figure{The moves through the front/back faces.}
\definexref{fig-moves4}{\the\figno}{fig}
\endinsert

\proof
As one would expect, this statement may be derived from (a refined 
version of) the transversality theorem. Details of the reasoning 
are very similar to those given by \'Swi\c atkowski 
in \cite{jacek}, with just a few changes.
The euclidean space $\r^3$ needs to be replaced with $S^3$,
and the plane $Q$---with the torus $T$. 
The possibility of a nontransversal
self-intersection of the diagram (cf.~\cite{jacek}, p.203)
is no longer excluded; the projection of a \lk\ on the torus
may have self-tangencies, provided that the condition of 
(IIb) 
is satisfied. Such a self-tangency yields the move of
\ref{fig-moves2b}.

Furthermore, in \'Swi\c atkowski's Definition~3.2 we 
require an additional transversality condition, namely
that the isotopy $H:S^1\times I \to S^3$ is transversal to
$c_f\cup c_b\subset S^3$. 
It implies that, if $H(S^1,t)\cap c_f\neq\emptyset$ for a given $t$,
then a move of type (IVa) takes place at $t$, and analogously
for $H(S^1,t)\cap c_b\neq\emptyset$ we have a move of type (IVb).
\qed

\bigskip
\rmk{Example}
The four \leg\ (un)knots depicted on \ref{fig-unknots} are \leg\ isotopic.

\midinsert
(a) $\Gamma_{1,0}$ with a zigzag; \hfill 
(b) the flying saucer; \hskip6mm \null
\medskip 
\centerline{\epsffile{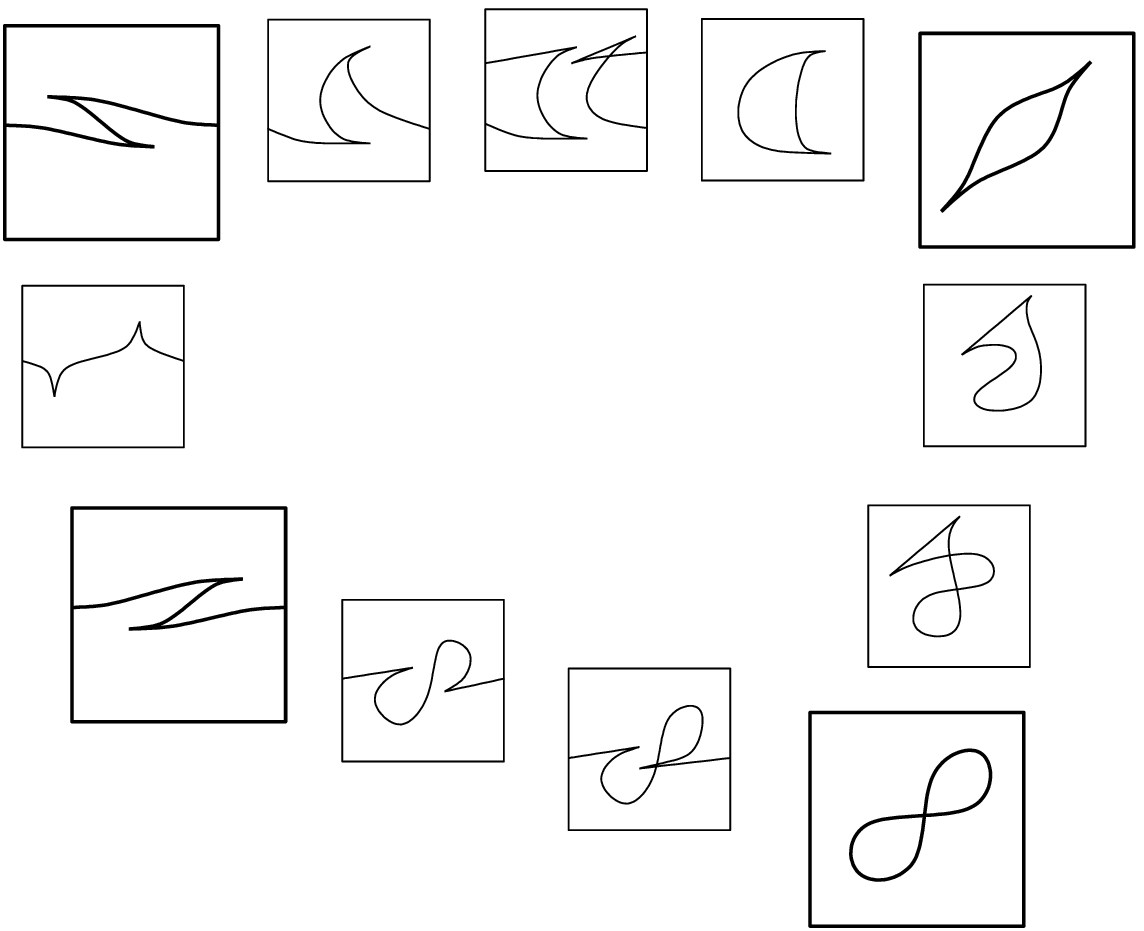}}
\vskip-2cm
{\parindent=0pt
\obeylines
(c) a \leg\ unknot 
in a small \nbd 
of the front face (caution:
it is {\it not} $\Gamma_{1,0}$ with the other zigzag);} 
\vskip4mm
\rightline{(d) a \leg\ unknot with cuspless front.}
\figure{Four isotopic \leg\ unknots.}
\definexref{fig-ufo}{\the\figno(b)}{fig}
\definexref{fig-unknots}{\the\figno}{figure}
\endinsert

\rmk{Example}
The two knots depicted at the top of \ref{fig-chekexples} are \leg\ isotopic.
Note that these knots considered as \lk s in the standard tight \cstr\
(known as {\it the Chekanov examples}) are not \leg\ isotopic.
\definexref{ex-cheka}{\here}{exple}

\midinsert
\centerline{\epsffile{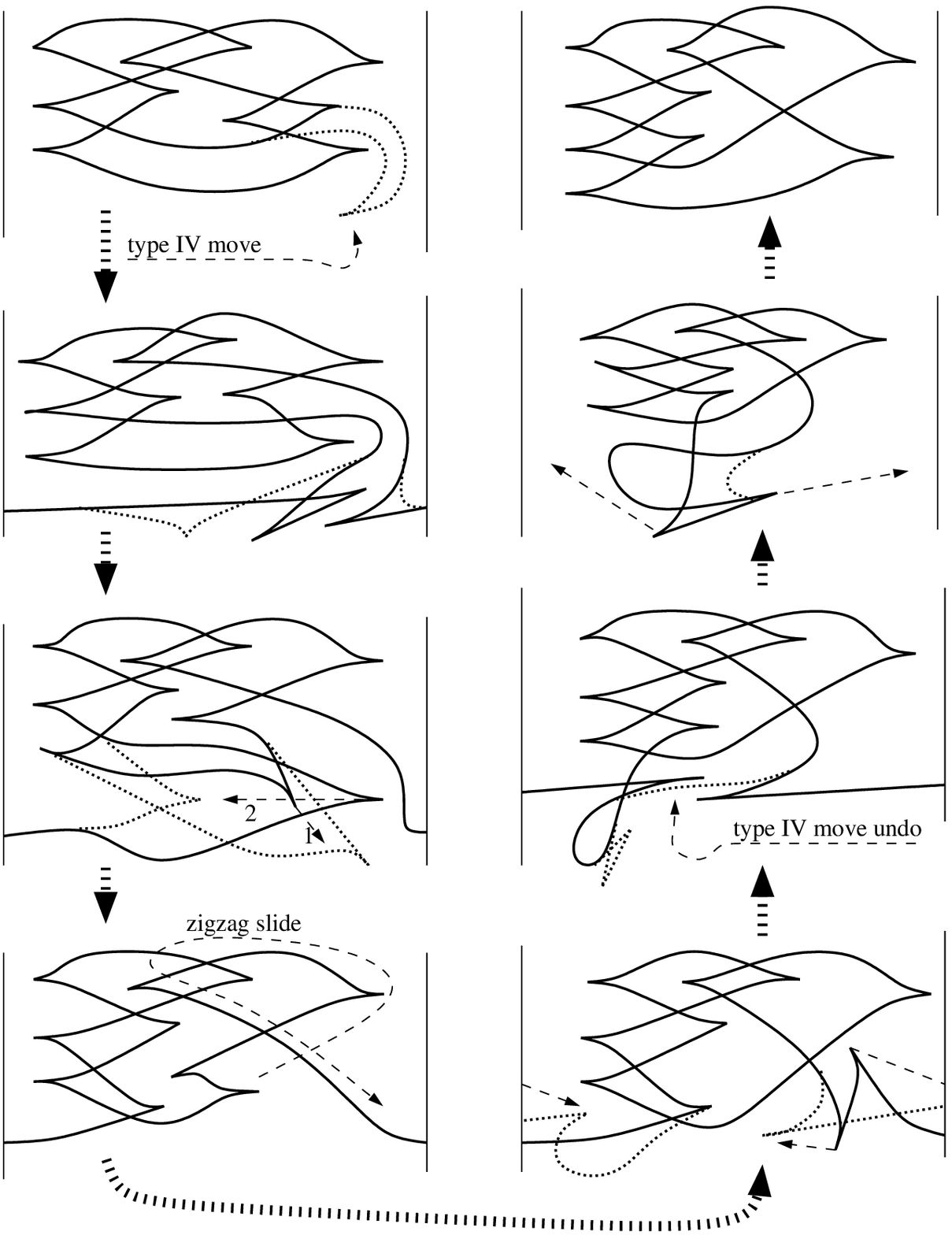}}
\figure{The Chekanov examples regarded as \lk s} 
\centerline{in an \ocstr\ are isotopic.}
\definexref{fig-chekexples}{\the\figno}{figure}
\endinsert

\bigskip
For a manifold $M$ with a co-oriented \cstr\ $\xi$
the contact embedding 
$$
j:(\r^3,\zeta_0)\hookrightarrow(M,\xi)
$$
is unique up to isotopy, because any such embedding can be
isotoped to one with the image in an arbitrarily small ball. 
When $M=S^3$ and $\xi=\zeta_1$, we get the front diagram of 
$j(K)$ simply transferring the diagram of $K$ to the torus
without any change of shape.
Of course, the image of a \leg\ isotopy is still 
a \leg\ isotopy; thus $j$ induces a map
$$j_*:\{\hbox{\it \lk\ types in }\zeta_0\}
\to\{\hbox{\it \lk\ types in }\xi\}.$$

\ref{ex-cheka} shows that this map does not have to be injective.
Neither is it necessarily surjective;
in the following we will call a \lk\ {\it standard}
if it lies in the image of $j_*$.

\sec{The Classical Invariants}
\definexref{sec-invts}{\the\secno}{section}

\bigskip
Now we recall the definitions of two 
classical (integer--valued) invariants of \leg\ isotopy
of knots, the {\it rotation} (or Maslov) number
and the {\it Thurston-Bennequin} invariant.
They were first introduced by Bennequin in \cite{bennequin}, 
end since then has been widely used for \lk s in the
standard tight \cstr\ on $\r^3$ (or $S^3$).
As we are mainly interested in \ocstr s, the known results
often do not apply directly to our situation and must be
reformulated and reproved. In particular, in \refs{propo-diagramrot}
and \refn{propo-diagramtb} (whose proofs constitute 
more than half of this section) we provide a way 
to calculate the values of invariants of a \lk\ in an \ocstr\ 
from its \fott.

Please note that all the knots considered in this section 
are {\it oriented}.

\rmk{Definition} 
\definexref{def-invts}{\here}{defin}
Let $K$ be a \lk\ in a \cstr\ $\zeta$.

\noindent {\bf (i)} \definexref{def-tb}{\here(i)}{defin}
Assume $K$ is homologous to zero.
Denote by $K'$ be the transverse knot obtained by small shift of $K$ 
in the direction normal to the contact planes. The linking number 
$lk(K,K')$ (i.e.~the number of times, counted with signs, $K'$ 
intersects a Seifert surface of $K$) is called the {\it Thurston-Bennequin 
invariant} of $K$ and denoted by $tb(K)$. Note that $tb(K)$ is indeed
determined by $K$ alone.

\noindent {\bf (ii)} \definexref{def-rot}{\here(ii)}{defin}
Assume that $\zeta$ is parallelized. 
The index (number of revolutions) of the tangent vector to $K$ with respect 
to the parallelization is called the {\it rotation} of the knot $K$ and 
denoted by $rot(K)$.

\medskip 
The \tbi\ is a particular case of the self-linking number
of a (null-homologous) framed knot; two topologically isotopic
\lk s have equal \tbi s if and only if they are isotopic as framed knots.
Sometimes (e.g.~in \cite{chernov:vassiliev})
the isotopy class of framed knots is used as a natural
(though not so convenient) generalization of the \tbi\ to homologically
nontrivial knots. For a numerical version of this generalization
(``affine self-linking number'') see \cite{chernov:aslk}.

The parallelization of the \cstr\ in \ref{def-rot}
does not have to be neither unique nor defined on the
whole manifold: rotations
of all knots contained in a fixed subset of the manifold
can be calculated with respect to one fixed parallelization.
Note also that
the rotation is actually an invariant of \leg\ {\it curves}\/
rather than knots, as it remains unchanged under a homotopy
through \leg\ curves with possible self-intersections.

Furthermore, the rotation changes sign when the knot orientation
is reversed.
On the other hand, the \tbi\ of a \lk\ does not depend on its orientation:
reversing the orientation of $K$ changes also the orientation of
$K'$, so the linking number remains intact.

\smallskip
Equipped with those two invariants (in addition to the topological knot type,
which is obviously also an invariant of \leg\ isotopy) 
we may ask a natural question:
is this set of invariants complete?

\thm{Question}{Let ${\cal L}_M$ be the set of \leg\ isotopy classess
of knots in $M$, ${\cal K}_M$ the set of topological knot types.
Consider the map ${\cal T}:{\cal L}_M\to{\cal K}_M\times\z\times\z$
associating with a \lk\ $K$ its topological type, rotation and
\tbi.
\item{\bf (a)} Is ${\cal T}$ injective?
\item{\bf (b)} Is ${\cal T}$ surjective?}
\definexref{que-main}{\here(b)}{question}

\medskip
In general, the answer to both parts is ``no''.
\item{\bf (a)}
The are known examples of topologically isotopic \leg\ knots
(including some in the standard tight \cstr\ on $\r^3$)
which have equal invariants but are not \leg\ isotopic 
(\cite{chekanov}, \cite{fraser:example}, \cite{ng:computable};
the pair of knots from \ref{ex-cheka} has this property as well).
In \ref{sec-mainthm} we discuss this question for
\ocstr s.
\smallskip
\item{\bf (b)}
It is known that in the standard tight \cstr\ on $\r^3$
the values of invariants of any \lk\ $K$ are related via
the congruence $tb(K)+rot(K)\equiv 1\pmod 2$.
\ref{lem-futa} allows to generalize this result to 
homotopically trivial \lk s in an arbitrary contact manifold.

\medskip
Moreover, in several contact manifolds there are various restrictions 
on the range of the values of invariants, 
the best known\numberedfootnote{Other similar in spirit bounds have been found by 
Fuchs and Tabachnikov in \futa.} being the Bennequin 
inequality \be.

\thmnamed{Theorem}{Bennequin inequality}{For $K$
a \lk\ in the standard tight \cstr\ we have
$$
tb(K)+rot(K)\leq-\chi(K),
$$
where $\chi(K)$ is the maximal Euler characteristic of a Seifert surface 
of $K$.}
\definexref{benineq}{\here}{thm}

Fuchs and Tabachnikov in \futa\ provide a way of reading the values
of $rot$ and $tb$ from the diagram data.
We give an analogous result for $(S^3,\zeta_1)$, but first we need to
introduce some notation.

Consider the homotopy type of a \fott\ as a curve in $T^2$, 
i.e.~an element $(\alpha,\beta)\in\pi_1T^2\simeq{\bf Z}^2$.
We choose the generators so that $(1,0)$ is the horizontal
line oriented left-to-right, and $(0,1)$---the vertical line
oriented upwards. We will denote the homotopy type of the
diagram of $K$ by $\bigl(\alpha(K),\beta(K)\bigr)$.


In a front diagram in the standard tight \cstr\ there are 
two sorts of (unoriented) cusps: those pointing to the left
and those pointing to the right; and each of them may be passed
either upward or downward. Also, the branches (pieces between cusps) 
of an oriented front are either traversed left-to-right or
right-to-left (and these two types alternate).

The same can be said about the \fotts, with proper understanding
of ``up'', ``down'', ``left'' and ``right''.
Pick the coorientation of $\xi$ pointing up near the front face.
It induces a normal vector field $\nu$ along the knot diagram. 
A cusp is {\it left} if it points to the left after  rotating 
so that $\nu$ is directed upwards, and {\it right} otherwise.
The same convention applies to the orientation of branches.
A cusp (left or right) is {\it ascending} if it is traversed in the
direction of $\nu$ and {\it descending} otherwise.
\ref{fig-cusps} shows all kinds of cusps.

\midinsert
\centerline{\epsffile{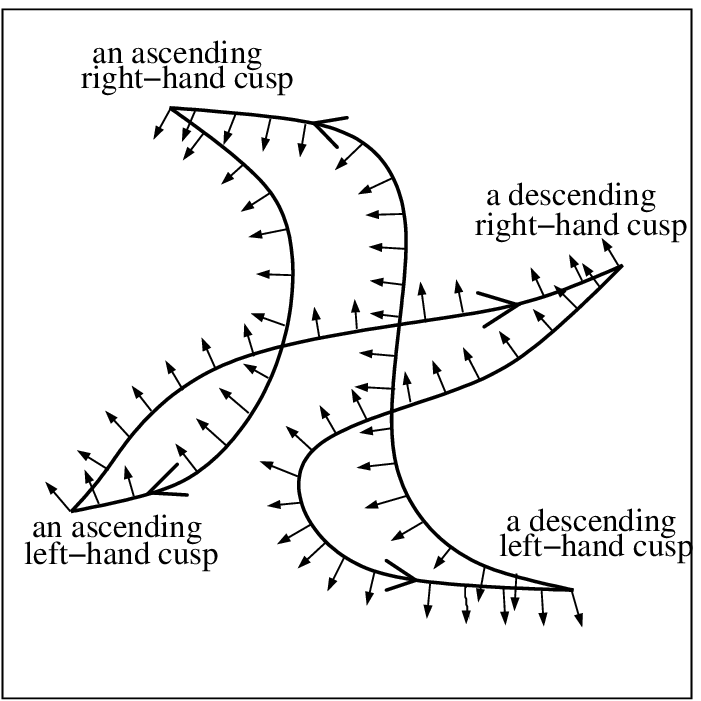}}
\figure{Types of cusps.}
\definexref{fig-cusps}{\the\figno}{figure}
\endinsert

\thm{Proposition}{Consider a \lk\ $K\subset(S^3,\xi)$
whose diagram $\pi(K)$ has the homotopy type 
$(\alpha(K),\beta(K))$ and contains $\cusp+(K)$ positive cusps
and $\cusp-(K)$ negative ones. Then}
$$
rot(K)={1\over2}\bigl(\cusp+(K)-\cusp-(K)\bigr)-\alpha(K)-\beta(K).
$$
\definexref{propo-diagramrot}{\here}{propo}

\proof
We need to pick a convenient parallelization of $\xi$. A nice candidate
would be a vector field which is constant in the $xyz$ coordinates of
\ref{cube}, say, pointing in the positive direction of the
$x$ axis (let us call this vector field $\cal X$). 
Unfortunately, $\cal X$ does not induce a vector field on the sphere.
We can, however, find a parallelization of $\xi$ on $S^3$
whose pull-back to the cube coincides with $\cal X$ in 
a region ${\cal C}_0\supset [\eps,1-\eps]^3$.
We may assume that the knot lies entirely in ${\cal C}_0$
except for a finite number of small segments at the points 
of intersection with the side faces.

Now we calculate the number of revolutions made by 
the tangent vector $\kappa$ 
within the region ${\cal C}_0$. Observe that $\kappa$ coincides
with $\cal X$ at each left positive cusp and each right
negative cusp, and that at the left positive cusps it is 
locally turning in the positive direction and at the right
negative ones---in the negative direction.
Therefore, the contribution to $rot(K)$ from the parts of $K$
contained in ${\cal C}_0$ is 
$\cusp+^{left}(K)-\cusp-^{right}(K)$.
Since $K$ must have exactly as many left cusps as right ones, 
$\cusp+^{left}(K)-\cusp-^{right}(K)=
\frac1/2\bigl(\cusp+(K)-\cusp-(K)\bigr)$.

Finally, we have to account for the contribution to $rot(K)$
of the segments at the intersection with the side faces.
Considering a small loop around the image of the front (resp.~back)
face it is easy to see that, since the parallelization extends to
$S^3$, the tangent vector necessarily makes one turn (see \ref{fig-rota}).
When we traverse the loop in the direction of the chosen generator,
the turn is taken in the negative direction.
A continuity argument shows that the same is true for any \lk\ 
crossing the side face. Counting with signs,
$K$ intersects the side faces $\alpha(K)+\beta(K)$ times.
Hence the formula.
\qed

\midinsert
\centerline{\epsffile{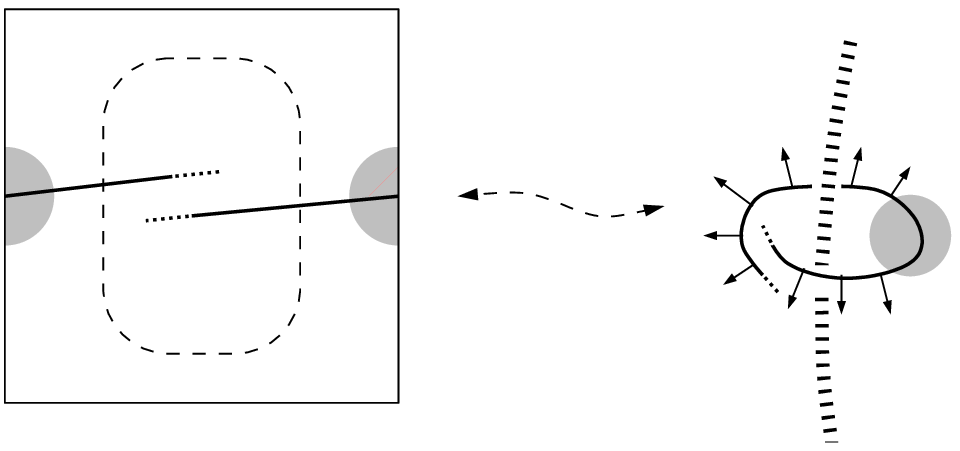}}
\figure{Crossing the side face contributes to rotation.}
\definexref{fig-rota}{\the\secno}{\the\figno}
\endinsert

\medskip
In order to formulate the results concerning the 
\tbi, we need to distinguish between positive and 
negative crossings. We adopt the usual convention, 
recalled in \ref{fig-crossings}.
However, the front diagrams of \lk s traditionally 
do not show the over/under crossing indication (since
it is redundant), which makes reading the sign of a 
crossing somewhat inconvenient. 

\midinsert
\centerline{\epsffile{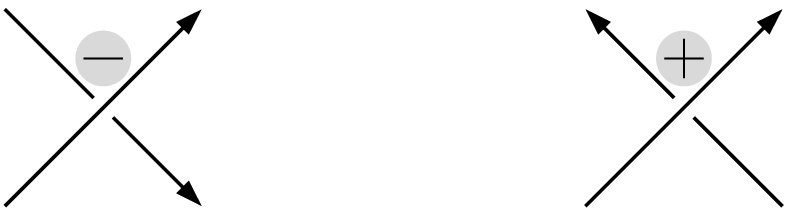}}
\figure{Positive and negative crossings.}
\definexref{fig-crossings}{\the\figno}{fig}
\endinsert

In the standard tight \cstr, a helpful observation 
is that crossings between a branch oriented left-to-right
and one oriented right-to-left is always negative, 
while a crossing between two branches of matching
orientations is positive. 
The same (with ``left'' and ``right'' understood as explained above)
is true for the \cstr\ $\xi$ and \fotts,
provided that the $x$ coordinates of the two branches
differ by less than $\pi$; if they differ by more than $\pi$,
the situation is opposite.

\medskip
The following lemma temporarily leaves the realm of contact topology.
For $K$ an arbitrary topological knot in 
$S^3\setminus(c_f\cup c_b)$
we may still consider its diagram (projection) on the torus. 
There will be no cusps anymore (not for a knot in general position),
but the homotopy type 
$\bigl(\alpha(K),\beta(K)\bigr)$, as well as the distinction between
positive and negative crossings, still make sense.

\thm{Lemma}{Suppose that 
\item{\rm (i)} $K_1$ and $K_2$ are topological knots in 
$S^3\setminus(c_f\cup c_b)$;
\item{\rm (ii)} the diagram $\pi(K_i)$ has the homotopy type
$\bigl(\alpha(K_i),\beta(K_i)\bigr)$;
\item{\rm (iii)} there are $\cross+(K_1,K_2)$ positive crossings 
and $\cross-(K_1,K_2)$ negative ones with one strand from 
$K_1$ and the other from $K_2$.
\hfill\break
\noindent Then the linking number (in the sphere) is given by the formula
$$lk(K_1,K_2)={1\over 2}\bigl(\cross+(K_1,K_2)-\cross-(K_1,K_2)
-\alpha(K_1)\beta(K_2)-\alpha(K_2)\beta(K_1)\bigr)
$$}
\definexref{lem-lk}{\here}{lem}

\proof
The formula 
$lk(K,L)={1\over 2}\bigl(\cross+(K_1,K_2)-\cross-(K_1,K_2)\bigr)$
holds for knots in $\r^3$ and their diagrams on the plane.
We can consider a knot in the 3-sphere 
as a knot in $\r^3=S^3\setminus\{*\}$.
To get its diagram on the plane from the one on the torus
it is enough to close up the loose ends which arise when we cut the
torus and unfold it to a square (see \ref{fig-closeup}).
Notice that all the lines connecting the horizontal loose ends overcross
those attached to the vertical ends. Hence the additional lines
contribute 
$\alpha(K_1)\beta(K_2)+\alpha(K_2)\beta(K_1)$ 
crossings counted with signs. 

The formula follows.
\qed

\midinsert
\centerline{\epsffile{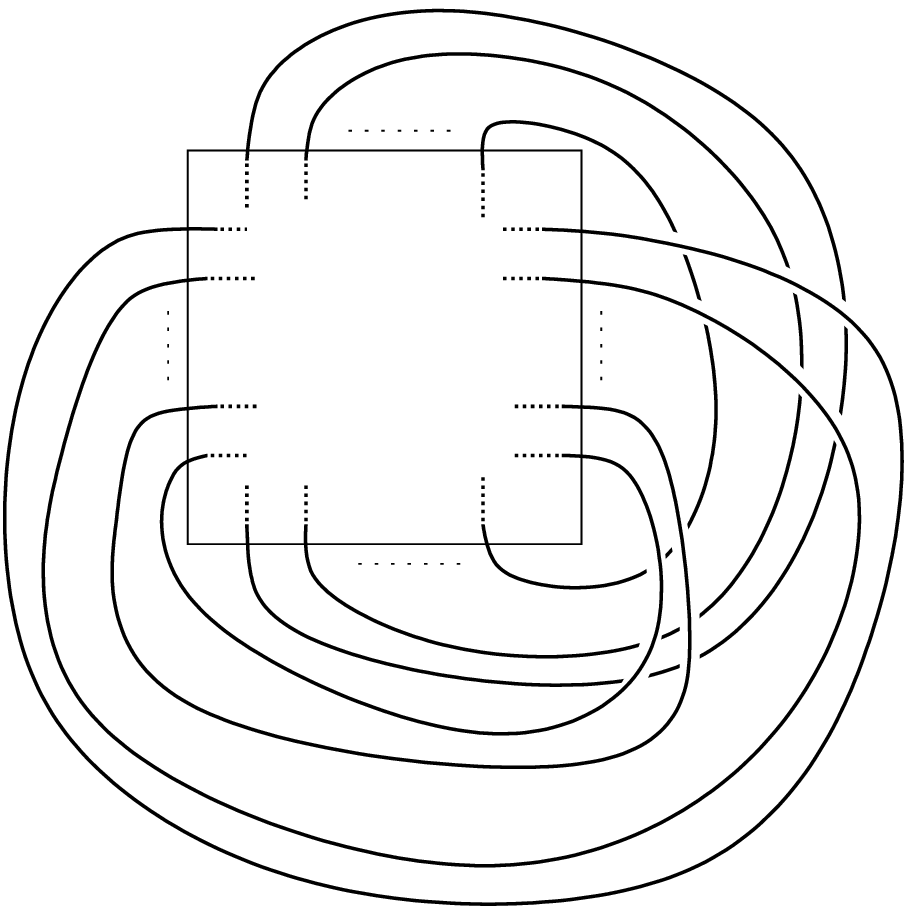}}
\figure{Closing up the diagram.}
\definexref{fig-closeup}{\the\figno}{fig}
\endinsert

\thm{Proposition}{Let $K$ be a \lk\ in $(S^3,\xi)$ whose diagram 
has the homotopy type $\bigl(\alpha(K),\beta(K)\bigr)$,
contains $cusp(K)$ cusps, $\cross+(K)$ positive crossings
and $\cross-(K)$ negative crossings. Then
$$tb(K)=\cross+(K)-\cross-(K)-{1\over2}\,cusp(K)-\alpha(K)\beta(K)
$$}
\definexref{propo-diagramtb}{\here}{propo}

\proof
By \ref{def-tb} and \ref{lem-lk}, 
$tb(K)=lk(K,K')=
{1\over 2}\bigl(\cross+(K,K')-\cross-(K,K')
-\alpha(K)\beta(K')-\alpha(K')\beta(K)\bigr)$, 
where $K'$ is the knot
obtained from $K$ by pushing it slightly in the direction of
the coorientation of $\xi$. 
The diagram (\fott) of $K'$
is also obtained by a small shift from that of $K$.
Therefore 
\item{$\bullet$} each positive crossing in the diagram of $K$
gives rise to two positive crossings between the diagrams
of $K$ and $K'$;
\item{$\bullet$} each negative crossing in the diagram of $K$
gives rise to two negative crossings;
\item{$\bullet$} each cusp in the diagram of $K$
gives rise to a negative crossing;

\noindent
and these are all the crossings between the diagrams of $K$ and $K'$.

Certainly, $\alpha(K')=\alpha(K)$ and $\beta(K')=\beta(K)$.

Thus,
$$\eqalign{lk(K,K')&=
{1\over 2}\,\bigl(\cross+(K,K')-\cross-(K,K')
-\alpha(K)\beta(K')-\alpha(K')\beta(K)\bigr)=\cr
&={1\over 2}\,\bigl(2\,\cross+(K)-2\,\cross-(K)-cusp(K)
-2\,\alpha(K)\beta(K)\bigr).}$$ 
The formula follows.
\qed

\bigskip

\refs{propo-diagramrot} and \refn{propo-diagramtb}
imply that the values of invariants are always uniquely determined
by the diagram of a \lk, even if the knot itself is not.
(This is the case, for example, with some of the knots $\g mn$
defined in \ref{ex-gammas}.) We want to calculate $rot$ and $tb$ of
$\g mn$'s; to this end it is necessary to choose an orientation
for each of them. We do it in an obvious way, requiring that 
$\alpha(\g mn)=m$ and $\beta(\g mn)=n$. 
Thus, $\g{-m}{-n}=\overline\g mn$ (i.e.~$\g mn$ with the
opposite orientation).

\rmk{Example}
With the orientation as described above, we have
$rot(\g mn)=-m-n$ and $tb(\g mn)=-mn$.

\medskip
Note that topologically any $\g mn$ is a torus knot
of type $(m,n)$; in particular, for $|m|=1$ (or $|n|=1$)
it is an unknot. Therefore for each integer $n$ we get
a \leg\ unknot $\g 1n$ of rotation number $-n-1$ and
\tbi\ $-n$. 
\bigbreak

\sec{The Classification Theorem}
\definexref{sec-mainthm}{\the\secno}{section}

In this section we formulate and prove the theorem 
announced in the introduction
as a generalization (with completely different methods of proof) 
of the main theorem
of \dy.

\thm{Theorem}{Let $\xi$ be an \ocstr\ on a 3-manifold $M$,
$\Delta$ an \ot\ disk in $(M,\xi)$, 
$K$ and $L$ \lk s in $M\setminus\Delta$.
Assume that $\xi$ is trivializable as a plane bundle.
Let $\cal K$ be the topological knot type of $K$.
Assume that:
\item{\rm (a)} there are infinitely many non-isotopic framed knots
in $\cal K$;
\item{\rm (b)} $K$ and $L$ are isotopic as framed knots
(with \leg\ framing); and 
\item{\rm (c)} $rot(K)=rot(L)$, where $rot$ is calculated
with respect to some parallelization of $\xi$.

\noindent Then $K$ and $L$ are \leg\ isotopic.}
\definexref{thm-main}{\here}{thm}

\rmk{Remark}
The assumption (a) is satisfied for all knots in a manifold
$M$ (not necessarily closed) which is not realizable as a \consum\ 
$M'\#(S^1\times S^2)$ (see Theorem 2.0.5 in \cite{chernov:aslk}).

For homologically trivial knots the assumption (b)
is equivalent to the conjunction of the following
two conditions:
\item{\rm (b1)} $K$ and $L$ are isotopic as topological knots, and
\item{\rm (b2)} $tb(K)=tb(L)$.
\bigskip

The proof of this theorem relies of two lemmas.
\ref{lem-makezigzag} explains the role of the \ot\ disks.
\ref{lem-futa} generalizes (to an arbitrary \cstr)
Theorem~4.4 of \futa; its main ingredient is
the notion of adding a zigzag, a.k.a.~stabilization,
also introduced by Fuchs and Tabachnikov in \futa.
They understand adding a zigzag as a special case of 
the operation of \consum mation of \lk s. For us, this approach 
would have serious disadvantages: the \consum\ of \lk s
in an arbitrary \cstr\ is not always well defined.
We will, therefore, content ourselves with a na\"\i ve diagrammatic
definition of adding a zigzag.\numberedfootnote{There 
is another fairly general notion which contains adding a
zigzag as a special case, namely the (\leg) satellite construction
(see \cite{ng:satellites}).}

\rmk{Definition} \definexref{def-stab}{\here}{defin}
Let $(M,\zeta)$ be a contact manifold, $K$ a \lk\ in it.
Pick a point $x_0\in K$; there is a \nbd\ $V\owns x_0$ such that 
$\zeta|_V$ is the standard \cstr\ on $\r^3$\numberedfootnote{A
{\it Darboux chart} in the nomenclature of \cite{chernov:vassiliev}.} 
(and $K\cap V$ a long \lk\ in it), thus it makes sense to talk
about the front diagram of $K\cap V$.
Now replace a small cuspless segment of $K\cap V$ 
with one containing $p$ pairs of ascending cusps 
and $q$ pairs of descending cusps (see \ref{fig-stab}).
We call the resulting knot {\it the $p,q$-stabilization}\/
of $K$ and denote it by $Z_{p,q}(K)$.

\midinsert
\centerline{\epsffile{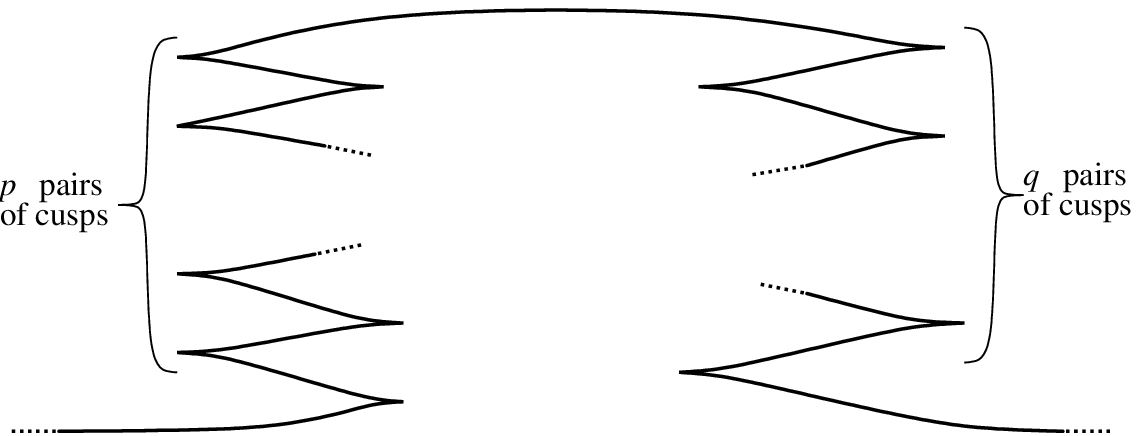}}
\figure{Adding a zigzag.}
\definexref{fig-stab}{\the\figno}{fig}
\endinsert

\rmk{Example} 
\ref{fig-zigzags} shows a few of
the \lk s obtained by adding zigzags to 
the standard ``flying saucer'' Legendrian unknot
(depicted in \ref{fig-ufo}). 
We will not distinguish between the operation of
adding a zigzag and the result of this operation
applied to the standard unknot, 
denoting these knots simply by $Z_{p,q}$.

\midinsert
\centerline{\epsffile{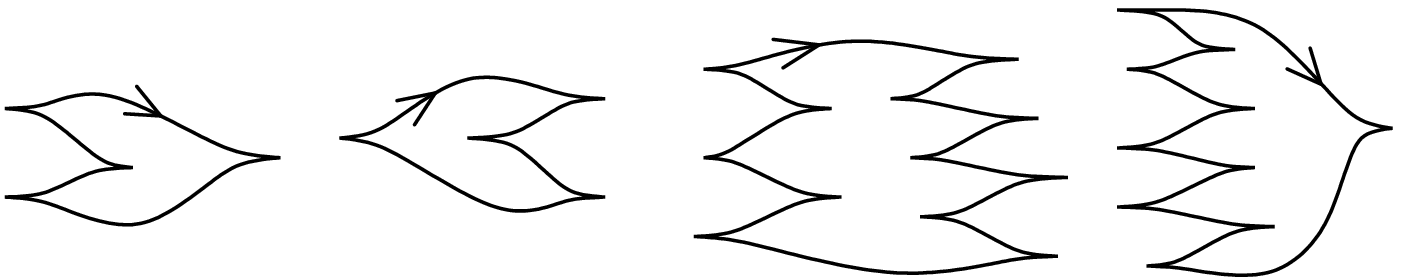}}
\noindent
\null\hskip1cm $Z_{2,1}$\hfill $Z_{1,2}$
\hfill $Z_{3,4}$ \hfill $Z_{5,1}$\hskip1cm\null
\figure{Examples of zigzags.}
\definexref{fig-zigzags}{\the\figno}{fig}
\endinsert

\medskip

One needs to show that $Z_{p,q}(K)$ is well defined and does not 
depend on the choice of the point $x_0$ where the new cusps are placed.
The classical argument, based on sliding pair of cusps along the knot
(cf.~\futa, Lemma~4.3 with Figures~17 and 18, which are recalled here 
as \ref{fig-cuspslide}), applies directly.


\midinsert
\centerline{\epsffile{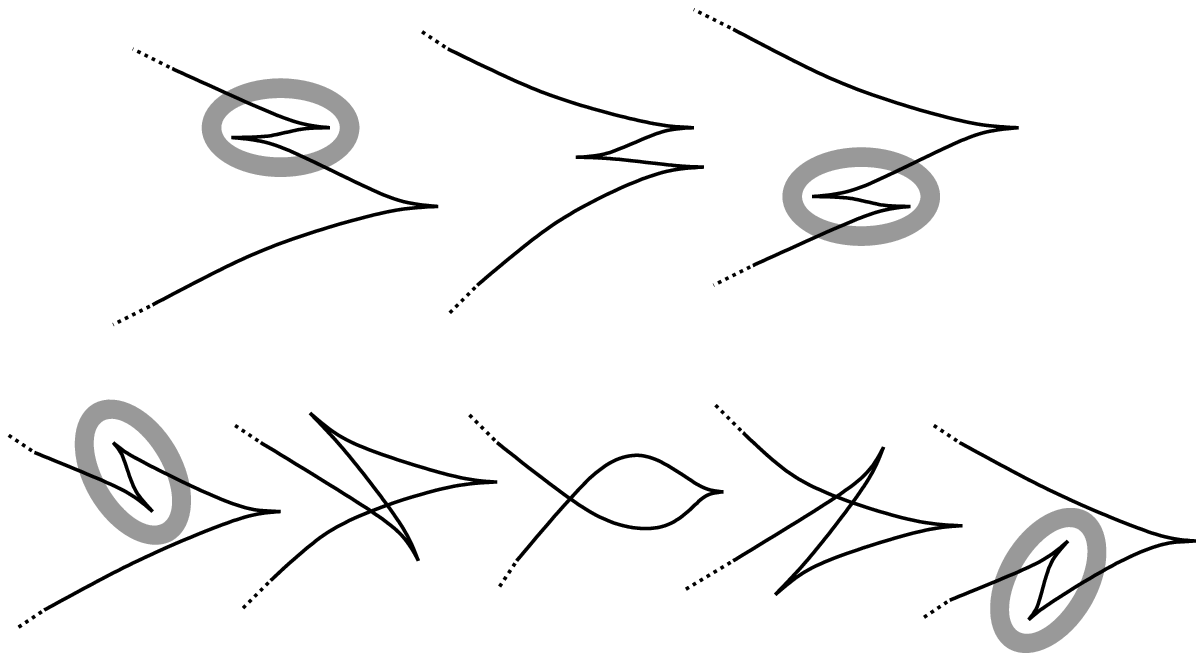}}
\figure{Sliding pair of cusps (elementary zigzags).}
\definexref{fig-cuspslide}{\the\figno}{fig}
\endinsert

The following proposition collects 
a couple of (obvious, but useful)
facts about adding zigzags.

\thm{Proposition}{\item{\bf (a)} 
$Z_{p,q}\bigl(Z_{p',q'}(K)\bigr)$ 
{\it is \leg\ isotopic to} $Z_{p+p',q+q'}(K)$;
\item{\bf (b)} 
$rot\bigl(Z_{p,q}(K)\bigr)=rot(K)+p-q$; and 
\item{\bf (c)}
the \leg\ framing of $Z_{p,q}K$ is the \leg\ framing of $K$ with $p+q$
negative twists added (for zero-homologous $K$ this means 
$tb\bigl(Z_{p,q}(K)\bigr)=tb(K)-p-q$.)} \qed
\definexref{prop-zig}{\here}{propo}

\thm{Lemma}{If $K$ and $L$ are two (long or short) \lk s 
in an arbitrary contact manifold $(M,\zeta)$ 
which are topologically isotopic,
then there exist two \lk s $K_Z=Z_{p,q}(K)$ and $L_Z=Z_{p',q'}(L)$ 
such that
$K_Z$ is \leg\ isotopic to $L_Z$.}
\definexref{lem-futa}{\here}{lemma}

\proof
We begin by choosing a topological isotopy between $K$ and $L$ 
which makes changes piece by piece only.
Namely, 
for $\{K_t, t\in [0,1]\}$, $K_0=K$, $K_1=L$,
we assume 
that there are finitely many points $t_0=0<t_1<t_2<\ldots<t_n=1$
and balls $B_1, B_2,\ldots, B_n\subset M$ such that
\item{$\bullet$} $\zeta|_{\textstyle{B_i}}$ is the standard tight \cstr\
on $\r^3$; 
\item{$\bullet$} for $t_{i-1}\leq t \leq t_i$, $K_t$ is constant outside $B_i$.

Now it is enough to study the isotopies
$\{K_t|_{\textstyle{B_i}}, t\in[t_{i-1},t_i]\}$ 
as isotopies of \lk s in the standard tight \cstr\ on $\r^3$.
By Theorem~4.4 of \futa,  $Z_{p_{i-1},q_{i-1}}(K_{t_{i-1}})$
is \leg\ isotopic to $Z_{p_i,q_i}(K_{t_i}\#Z_{p_i,q_i})$ for some
$p_{i-1}, q_{i-1}, p_i, q_i$.
An easy inductive argument ends the proof.
\qed

\thmnamed{Lemma}{creating zigzags}{Let $V(\Delta)$ be a small \nbd\ 
of an \ot\ disk $\Delta$ in $M$, 
$L$ a \lk\ in $M\setminus \Delta$.
Assume that (under the identification of $V(\Delta)$ with a subset
of $(S^3,\zeta_1)$ described above) the projection of the curve
$\ell=L\cap V(\Delta)$ onto the torus is a straight line segment
(as in the top part of \ref{fig-makezigzag}, where the shaded area
symbolizes $V(\Delta)$.\numberedfootnote{Note that, even though 
the \fott\ of the top left part of \ref{fig-makezigzag} 
corresponds to two different \leg\ curves, 
the condition that $L$ is in $M\setminus \Delta$
actually excludes one of them.}
Then for any $p$ and $q$
there is a \lk\ $\ell_{p,q}$ in $U$ coinciding with $\ell$
outside a compact set and such that
$Z_{p,q}(K\setminus\ell)\cup\ell_{p,q}$ is \leg\
isotopic to $K$.} 
\definexref{lem-makezigzag}{\here}{lemma}

\proof
We begin by checking this statement for elementary zigzags. 
\ref{fig-makezigzag} shows an isotopy between $L$ and a \lk\
coinciding outside $V(\Delta)$ with $Z_{0,1}(L)$ or $Z_{0,1}(L)$, 
depending on the orientation of $\ell$.

In order to produce the other elementary zigzag, 
first perform the type I Reidemeister move on $\ell$,
so that it crosses some smaller \nbd\ of $\Delta$
three times: twice in the same direction as before and once
(between cusps) in the opposite direction.
By a small isotopy we can ensure that the slope of the segment
between cusps is the smallest, see \ref{fig-uturn}.
Then the corresponding segment of the curve lies closest to the disk:
we can replace $V(\Delta)$ with a smaller \nbd\ $V'(\Delta)$ which contains
only the segment going in the opposite direction to that of the 
original $\ell$.

\midinsert
\centerline{\epsffile{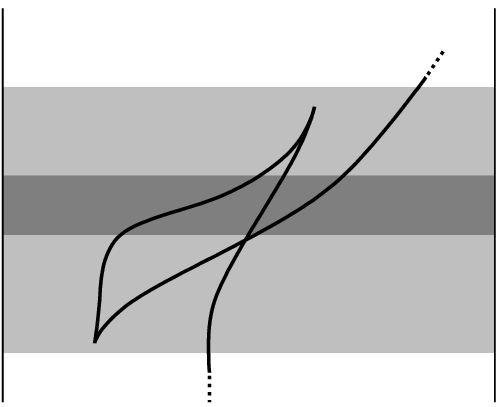}}
\figure{Changing the direction of $\ell$}
\definexref{fig-uturn}{\the\figno}{figure}
\endinsert

In order to prove the lemma for an arbitrary zigzag $Z_{p,q}$,
use $p+q$ disjoint slices of $V(\Delta)$ (\nbd s of $p+q$ parallel copies
of $\Delta$), to produce $p$ copies of $Z_{1,0}$ and $q$ copies
of $Z_{0,1}$. Then slide the newly created zigzags along the knot
until they leave the whole $V(\Delta)$.
\qed

\midinsert
\centerline{\epsffile{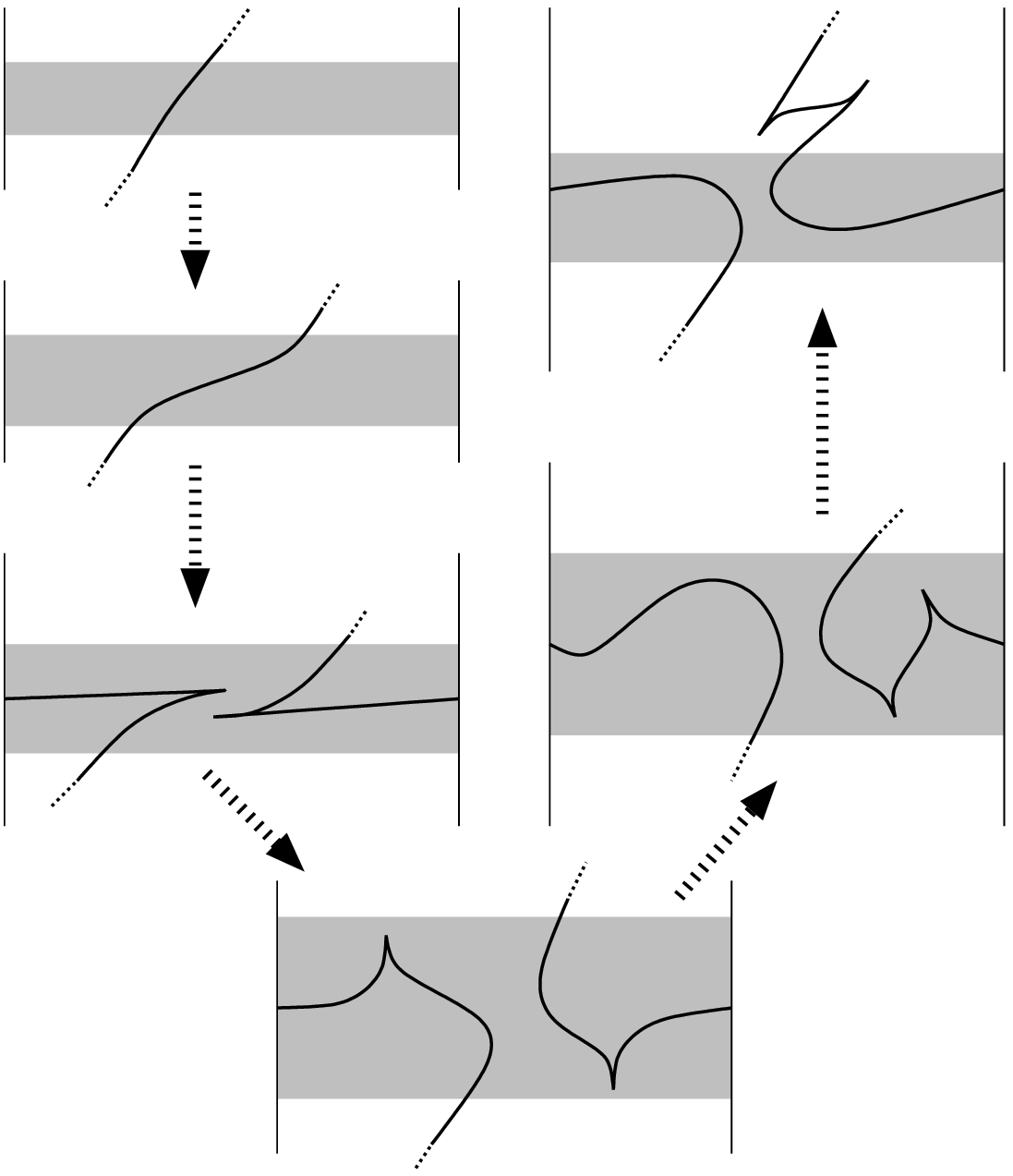}}
\figure{Making a zigzag (out of nothing).}
\definexref{fig-makezigzag}{\the\figno}{fig}
\endinsert

\bigskip
\noindent
{\it Proof of \ref{thm-main}:}
There is a \leg\ isotopy moving $K$ in $M\setminus \Delta$ 
so that it coincides with $L$ along a small segment, 
which we will call $\ell$.
By assumption, there is a topological isotopy between $K$ and $L$;
without loss of generality we may assume that the isotopy is constant 
on $\ell$. 

Pick a \nbd\ $V$ containing $\Delta$ and such that 
$\overline V\cap(K\cup L)=\ell$.
By \ref{lem-futa}, 
$K\setminus\ell$ and $L\setminus\ell$ 
(considered as {\it long} \lk s in $M\setminus V$)
are \leg\ isotopic modulo adding zigzags, 
i.e.~$Z_{p_K,q_K}(K\setminus\ell)$ and 
$Z_{p_L,q_L}(L\setminus\ell)$ 
are \leg\ isotopic, thus also
$Z_{p_K,q_K}K\simeq Z_{p_L,q_L}L$
as \lk s in $M$.

By \ref{prop-zig}(b) we have then
$rot(K)+p_K-q_K=rot(L)+p_L-q_L$ 
By assumption \refn{thm-main}(c) $rot(K)=rot(L)$,
thus $$p_K-q_K=p_L-q_L \eqno{(*)}$$.

Adding a $p,q$-zigzag introduces $p+q$ extra
negative twists to the \leg\ framing.

The \leg\ isotopy between $Z_{p_K,q_K}K$ and $Z_{p_L,q_L}L$
induces an isotopy of framed knots.
Since \leg\ framings of $K$ and $L$ also coincide (see \refn{thm-main}b),
we see that the effect of adding $p_K+q_K$ negative twists
to the \leg\ framing of $K$
must be the same as that of adding $p_L+q_L$ ones.
This implies that $$p_K+q_K=p_L+q_L,\eqno{(**)}$$ 
provided that the condition \refn{thm-main}(a) is satisfied.

From (*) and (**) it follows that 
$p_K=p_L$ and $q_K=q_L$; for simplicity,
we will denote these numbers by $p$ and $q$ respectively.

By \ref{lem-makezigzag},
$K\simeq Z_{p,q}\big(K\setminus\ell)\cup\ell_{p,q}$ 
and
$L\simeq Z_{p,q}(L\setminus\ell)\cup\ell_{p,q}$;
but we have also
$Z_{p,q}(K\setminus\ell)\simeq Z_{p,q}(L\setminus\ell)$.
This proves $K\simeq L$.
\qed

\medskip

\rmk{Definition}
A \lk\ is called {\it loose} if its complement is \ot.

\smallskip
For example, all standard knots (as defined at the end of 
\ref{sec-lks}) are necessarily loose; also any knot whose front on the torus 
does not cross a fixed meridian of the torus (as on \ref{fig-afront})
is loose. 

On the other hand, in \dy\ the author proved that the knot $\Gamma_{-1,1}$
is non-loose. It have been believed a sort of one-time phenomenon
(even the name originally proposed by Eliashberg and Fraser in 
\cite{eliashberg-fraser} was ``exceptional knots'').
Etnyre and Ng conjecture that $\Gamma_{-1,1}$ (is isotopic to 
$\Gamma_{-1,1}$ and) is the only non-loose \leg\ unknot
that can be found in any \ocstr\ on $S^3$ (Conjectures 41 and 42
in \cite{problems}).
This is not true; in \ref{prop-nonloosegammas} we prove that 
all the knots $\Gamma_{m,n}$ for which $mn<0$ are non-loose.
In particular, $\Gamma_{-1,n}$ (for $n=1,2,3\ldots$) are
all different (distinguishable by $rot$ and $tb$) non-loose
\leg\ unknots.
  
The above \ref{lem-makezigzag} (thus also \ref{thm-main}) deals with 
loose \lk s only. This assumption is, of course, essential:
two topologically isotopic \lk s with the same values of rotation
and \tbi\ may be \leg\ non-isotopic simply because one of them
is loose while the other is not (cf.~\ref{ex-kozuch}).
But even for two loose knots \ref{thm-main} is not necessarily
applicable. Namely, there exist two loose \lk s 
for which there is no \ot\ disk disjoint with {\it both}\/ of them.

\thm{Proposition}{Every
\ot\ disk in $(S^3,\zeta_1)$ intersects
either $\Gamma_{0,1}$ or $\Gamma_{1,0}$.}
\definexref{propo-nonloose}{\here}{propo}

\proof
We want to show that the open manifold 
$M=S^3\setminus(\Gamma_{0,1}\cup\Gamma_{1,0})$ is tight;
it is enough to prove tightness of its universal cover $\widetilde M$
(because a hypothetical \ot\ disk in $M$ would lift to $\widetilde M$).
The idea of this proof is the same as employed in \dy\
in the proof of Lemma 4.2.1: we will find a contactomorphism of 
$\widetilde M$ with an open subset of $\r^3$ with the tight 
\cstr\ $\zeta'=\ker(\cos x\,dz-\sin x\,dy)$ described in \ref{ex-stdprim}
(a subset of a tight manifold is, of course, also tight).

\medskip
\noindent
{\bf Step 1:} {\it Topological description}
\smallskip

We will see $M$ as made of three parts: two solid tori
$M_f=M\cap\pr_{\cal C}\bigl([0,\frac1/3]\times[0,1]^2\bigr)$ and 
$M_b=M\cap\pr_{\cal C}\bigl([\frac2/3,1]\times[0,1]^2\bigr)$, and 
the ``middle slice'' $M_m=M\cap\pr_{\cal C}\bigl([\frac1/3,\frac2/3]\times[0,1]^2\bigr)$,
glued together along two annuli: $A_0=M_f\cap M_m$ and $A_1=M_b\cap M_m$.
The middle slice is diffeomorphic to $T^2\times[0,1]$ with a meridian
removed from one component of the boundary and a latitude---from the other;
hence its universal cover $\widetilde{M_m}=\r^2\times[0,1]\setminus L$,
where $L$ consists of two infinite families of parallel lines, one in the
universal cover of each boundary component of $T^2\times[0,1]$.
In $\widetilde{M_m}$ the annulus $A_0$ (resp.~$A_1$) lifts to the union
of infinitely many stripes. Along each stripe we need to attach a copy of
$\widetilde{M_f}$ (resp.~$\widetilde{M_b}$), which is diffeomorphic
to $H\times\r$ for $H$ a half-disk.\numberedfootnote{From the description
of $M$ we see that $H$ is $D^2$ with a boundary point removed;
we choose to see it as a half-disk to make the gluing picture easier
to understand.}
Thus, we get the space pictured in \ref{fig-walki}.

\midinsert
\centerline{\epsffile{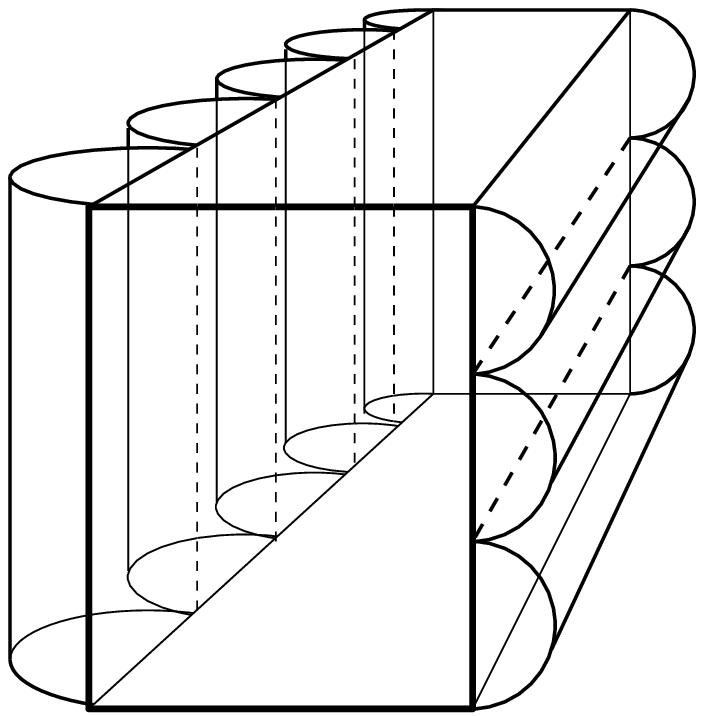}}
\figure{The universal cover $\widetilde{M}$.}
\definexref{fig-walki}{\the\figno}{fig}
\endinsert

\medskip
\noindent
{\bf Step 2:} {\it Contact structure}
\smallskip

First, notice that $\widetilde{M_m}$ can be contact embedded 
into $(\r^3,\zeta')$ so that its image consists
of the slice $\{0\leq x\leq \frac{$\pi$}/2\}$ with the subset $L$
(equal to $\{x=0,y\in\z\}\cup\{x=\frac{$\pi$}/2,z\in\z\}$)
removed.

Now we need to recall the construction used in the proof of 
Lemma 4.2.1 in \dy. We considered the universal cover of 
$S^3\setminus\Gamma_{-1,1}$ (with the same \cstr\ $\zeta_1$),
which again could be described as the union $\widetilde{M_0}\cup\widetilde{M_1}$,
with $M_0=\bigl(S^3\setminus\Gamma_{-1,1}\bigr)
\cap\pr_{\cal C}\bigl([0,\frac1/2]\times[0,1]^2\bigr)$ 
and  $M_1=\bigl(S^3\setminus\Gamma_{-1,1}\bigr)
\cap\pr_{\cal C}\bigl([\frac1/2,1]\times[0,1]^2\bigr)$,
glued together along an annulus $A$.
Then, we constructed a contact embedding $J:\widetilde M\hookrightarrow(\r^3,\zeta)$
(where $\zeta$ is the contact structure defined in \ref{ex-std}),
starting by setting its values on $\widetilde A$ so that $J(\widetilde A)=\{x=0\}$.
Thus $\widetilde{M_0}$ maps into the half-space $\{x\geq 0\}$.

Since the spaces $\widetilde{M_f}$ and $\widetilde{M_b}$ are both
contactomorphic to $\widetilde{M_0}$, we have contact embeddings
$j_f:(\widetilde{M_f},\zeta_1)\to(\r^3,\zeta)$ and 
$j_b:(\widetilde{M_b},\zeta_1)\to(\r^3,\zeta)$. This is not exactly
what we want:
\item{$\bullet$} the contact structure on the target $\r^3$ 
should be $\zeta'$ rather than $\zeta$;\numberedfootnote{Those two are,
of course, isomorphic; but since we are in the process of constructing
a specific embedding, we need proceed with greater exactness than
just up to a \cmm.} and
\item{$\bullet$} the boundary of $\widetilde{M_f}$ (resp.~$\widetilde{M_b}$)
should not be mapped onto the whole plane, but onto a single stripe.

\noindent
To get rid of both problems at once it is enough to notice that
each of the subsets $U_n=(\frac{$\pi$}/2,0]\times(n,n+1)\times\r$ and 
$V_n=[\frac{$\pi$}/2,\pi)\times\r\times(n,n+1)$ in $(\r^3,\zeta'_0)$
is contactomorphic to the half-space $\{x\geq 0\}$ in $(\r^3,\zeta)$
(and the \cmm, of course, maps the boundary of $U_n$ and $V_n$ 
onto the plane $\{x=0\}$).

Thus we get an embedding of a copy of $\widetilde{M_f}$ into each $U_n$
and a copy of $\widetilde{M_b}$ into each $V_n$. The union of boundaries
$\bigcup_n\partial U_n\cup\bigcup_n\partial V_n$ coincides with 
the boundary of the image of $\widetilde{M_m}$ via its previously constructed embedding
into $\{0\leq x\leq \frac{$\pi$}/2\}\subset (\r^3,\zeta'_0)$.
The flexibility which we enjoy when we start to construct the maps 
$j_f$ and $j_b$ is sufficient to ensure that the embedding of  
$\widetilde{M_m}$ and those of (the infinitely many copies of) 
$\widetilde{M_f}$ and $\widetilde{M_b}$ form together a contact embedding
$(\widetilde M,\zeta_1)\hookrightarrow(\r^3,\zeta'_0)$.
\qed

Thus we have left unanswered the following question.

\rmk{Question}
Are $\Gamma_{0,1}$ and $\Gamma_{1,0}$
\leg\ isotopic?

From \ref{thm-main} alone we can draw no conclusion in this matter;
in fact, the answer is unknown to the author.

\sec{The Catalog of \leg\ Knots}
\definexref{sec-zoo}{\the\secno}{section}

Let $\xi$ be an \ocstr\ on $S^3$. 
In this section we concentrate on existence results
for \lk s in $\xi$.
Our goal is to give
the answer (as comprehensive as possible) to \ref{que-main},
using the results of \refs{sec-lks} and \refn{sec-mainthm}.

\thm{Theorem}{For any topological knot type $\cal K$ 
and any pair of integers $r$ and $t$ such that $r+t \equiv 1 \pmod 2$,
there is a \lk\ $L$ in $(S^3,\xi)$ realizing $\cal K$ 
with $rot(L)=r$ and $tb(L)=tb$.}
\definexref{thm-epi}{\here}{thm}


\thm{Lemma}{For any pair of integers $r$ and $t$ such that 
$r+t \equiv 1 \pmod 2$ there is a \leg\ unknot $U=U(r,t)$
in $\xi$ with $rot(U)=r$ and $tb(U)=tb$.}
\definexref{lem-utr}{\here}{lemma}

\proof
We will draw front diagrams of the desired unknots.
Their appearance depends on the actual values of $r$ and $t$ as follows:

\smallskip
\item{(a)} For the zigzag $Z_{p,q}$ (see \ref{fig-zigzags}) 
we have $rot(Z_{p,q})=p-q$ and $tb(Z_{p,q})=-p-q-1$ (with the orientation
so chosen that there is $2p+1$ ascending and $2q+1$ descending cusps);
the pairs of integers $(r,t)$ that can be obtained this way 
are all those satisfying Bennequin inequality for unknots, 
i.e.~$t<0$ and $|r|\leq|t|$.

\smallskip
\item{(b)} The rotation of the knot whose front diagram on the cylinder can be seen
in the left part of \ref{fig-utr} is equal to $p-q$ and its \tbi\ is $p+q-1$;
therefore such knots realize all pairs $(r,t)$ with $t\geq 0$ and 
$|r|\leq|t|$;

\smallskip
\item{(c)} Finally, the knot whose front diagram is depicted in the 
right part of \ref{fig-utr} has the \tbi\ equal to $p-q$,
while its rotation may be $p+q$ or $-p-q$, depending on the orientation;
such knots exhaust all the remaining pairs $(r,t)$, i.e.~those where
$|r|\geq|t|$.

\midinsert
\centerline{\epsffile{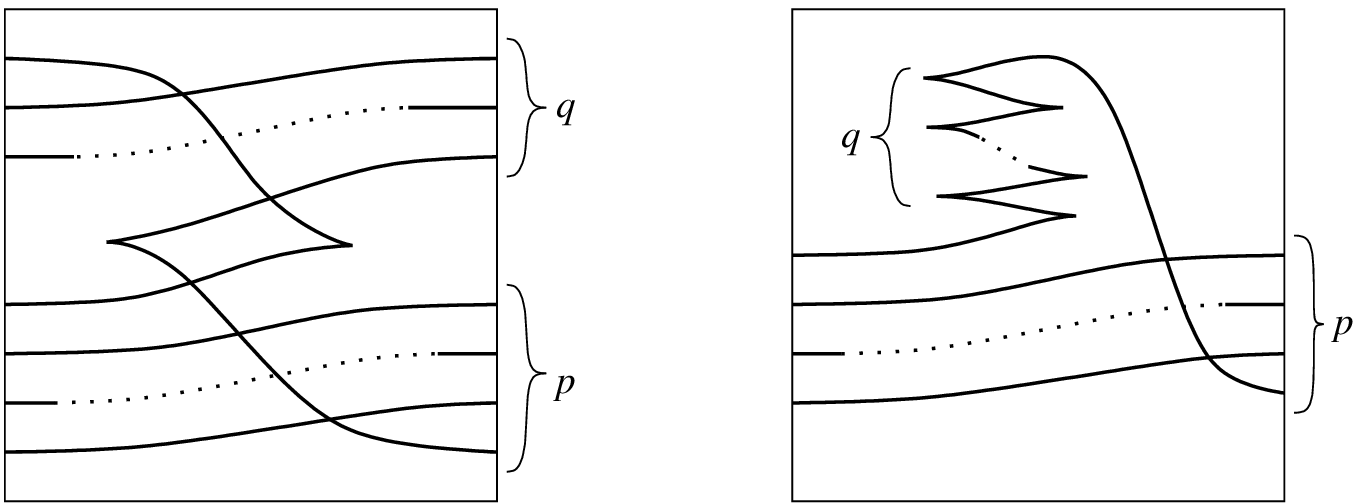}}
\figure{\leg\ unknots with $rot=0$ and $tb=t\geq 0$.}
\definexref{fig-utr}{\the\figno}{fig}
\endinsert

\bigskip
{\it Proof of \ref{thm-epi}:}

First, fix an \ot\ disk with its \nbd\ $V(\Delta)$ 
and choose a \lk\ $L_0$ realizing the topological knot type 
$\cal K$ so that the front diagram of its 
intersection with $V(\Delta)$ is a smooth arc as 
in the left part of \ref{fig-consum}
(the curve labeled $L_0$).
This is easily done, since any curve can be $C_0$-approximated
by a \leg\ one.

Next, take a \leg\ unknot $U$ with $rot(U)=r-rot(L_0)$ and
$tb(U)=t-tb(L_0)-1$, as constructed in the proof of \ref{lem-utr}.
We may assume that $U$ is contained in a smaller \nbd\ of
the \ot\ disk $V'(\Delta)\subset V(\Delta)$.

Then join $L_0\cap V(\Delta)$ with the lowest strand of $U$
in the way depicted in \ref{fig-consum} (in the case orientations
do not match, first perform the dovetail move on one of the curves).
We will denote the resulting knot by $L$.

\midinsert
\centerline{\epsffile{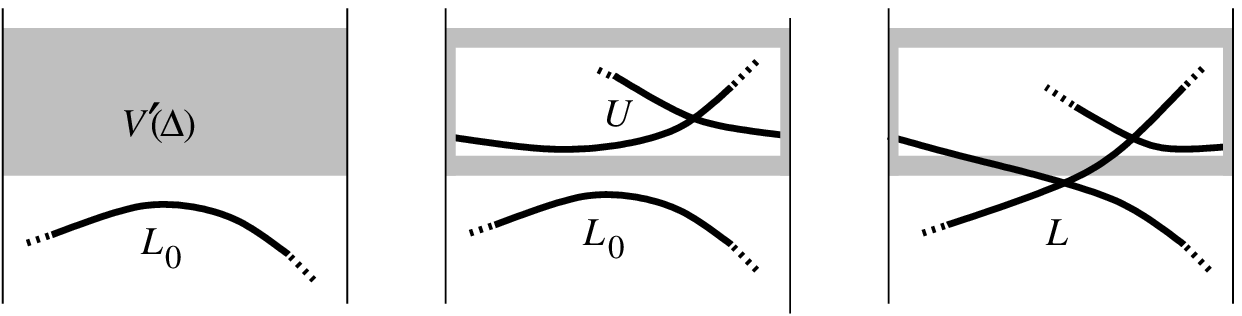}}
\figure{Constructing a \lk\ with pre-set invariants.} 
\definexref{fig-consum}{\the\figno}{figure}
\endinsert

The operation defined above resembles the diagrammatic definition
of the \consum mation of \lk s given by Etnyre and Honda in 
\cite{etnyre-honda:consums}.
The reason why we do not introduce it as such is that 
the \consum\ of \lk s may not be well-defined in general setting.
It retains, however, the following properties of the \consum mation 
(cf.~\cite{etnyre-honda:consums}, Lemma 3.3; note that 
it is proved without using the fact 
that the result of the \consum mation
is independent of the choice of the point where the operation is performed):
\item{(i)}
the topological knot type of $L$ is still $\cal K$, 
because $L$ is the {\it topological} \consum\
of $L_0$ with the unknot;
\item{(ii)}
$rot(L)=rot(L_0)+rot(U)=r$; and
\item{(ii)}
$tb(L)=tb(L_0)+tb(U)+1=t$.

Therefore $L$ satisfies all the conditions we require.
\qed

\medskip

All the \lk s constructed in the proof of \ref{thm-epi}
are loose.
We can, however, formulate some classification results 
specifically for non-loose knots.
\refs{prop-bensw} and \refn{prop-antybenneq} are corollaries
of the Bennequin inequality.
The first one of them has been pointed out to the author
by Jacek \hbox{\'Swi\c at}kowski.

\thm{Proposition}{Let $\xi$ be an \ocstr\ on a manifold $M$,
$K$ a non-loose \lk\ in $\xi$. Then
$-|tb(K)|+|rot(K)|\leq -\chi(K)$.}
\definexref{prop-bensw}{\here}{propo}

\proof
We want to apply the Bennequin inequality to a \lk\ $K'$
in the complement of $K$ (which is a tight manifold, 
since $K$ is non-loose).\numberedfootnote{A similar observation 
(though never explicitly expressed)
must have prompted Etnyre and Ng to make the ``final remark''
in Section~6 in \cite{problems}:
{\sl any \lk\ that violates the Bennequin inequality is automatically 
loose}.
This statement is not correct; \ref{prop-bensw} allows for 
non-loose knots violating the Bennequin inequality
(in the $tb>0$ case).
\ref{prop-antybenneq} shows that this correction is essential.}

We obtain $K'$ modifying 
the knot $K_\varepsilon$
used in \ref{def-tb} to define the \tbi\ so that it becomes 
homologically trivial in the complement of $K$.
The modification will be performed locally,
in a small tight ball, where the front diagrams of
$K$ and $K_\varepsilon$ are two parallel line segments.
The further procedure depends on the sign of 
$tb(K)$ as follows.
\smallskip

\noindent $tb(K)\leq0$
\item{}Replace $K_\varepsilon$ by $K'$ as in \ref{fig-benswiso}.
(In particular, for $tb(K)=0$ we take $K'=K_\varepsilon$.)
Here $K'$ is in fact \leg\ isotopic to $K$,
therefore obviously $tb(K')=tb(K)$ and $rot(K')=rot(K)$.

\noindent $tb(K)>0$
\item{} In this case, 
let $K'$ be the knot depicted in \ref{fig-benswzigzag}.
By \refs{propo-diagramrot} and \refn{propo-diagramtb}, 
$tb(K')=tb(K)-2\,tb(K)=-tb(K)$
(since there are $4\,tb(K)$ newly added cusps)
and $rot(K')=rot(K)$ (since there are equal numbers
of ascending and descending cusps among the new ones).

\midinsert
\vskip1cm
\centerline{\epsffile{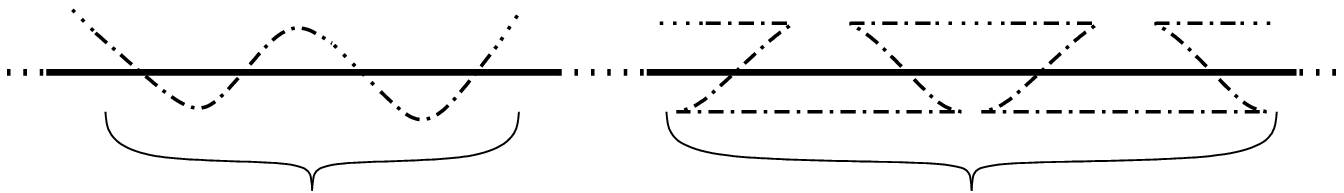}}
\hskip2cm $|tb(K)|$ crossings \hfill $tb(K)$ crossings \hskip2cm \null
\vskip7mm
\hskip2cm (a) $tb(K)\leq 0$ \hfill (b) $tb(K)>0$ \hskip2cm \null
\figure{Construction of the knot $K'$ in \ref{prop-bensw}.}
\definexref{fig-benswiso}{\the\figno(a)}{figure}
\definexref{fig-benswzigzag}{\the\figno(b)}{figure}
\vskip1cm
\endinsert

Thus in both cases we have $rot(K')=rot(K)$ and
$$
tb(K')=
\left\{ {tb(K)\atop -tb(K)} \quad {\hbox{ for }tb(K)\leq0\atop 
\hbox{ for }tb(K)>0} \right\}= -|tb(K)|.
$$
Moreover, since $\chi$ is a topological invariant, 
$\chi(K')=\chi(K)$.
Therefore, by \ref{benineq} for $K'$, 
$$
-|tb(K)|+|rot(K)|\leq -\chi(K).
$$

\thm{Proposition}{Let $\xi$ be an \ocstr\ on $S^3$, 
$K$ a non-loose \leg\ {\rm un}knot in $\xi$.
Then $tb(K)>0$ and $|rot(K)|<tb(K)$.}
\definexref{prop-antybenneq}{\here}{propo}

\proof
Suppose $K=\partial D$ ($D$ an embedded disk) is a \leg\ unknot
in $\xi$ and $tb(K)=0$. Then the knot $K'$ defined in the proof
of \ref{prop-bensw} is in fact isotopic to $K$.
On the other hand, $K'$ is a \leg\ unknot in $S^3\setminus D$.
The complement of the disk $D$ is homeomorphic to $\r^3$, and 
tight (because $K$ is non-loose).
Thus we can consider $K'$ as a \lk\ in the standard tight 
\cstr\ on $\r^3$.
Choose a tight ball in $S^3\setminus D$ so that there is an 
\ot\ disk disjoint with it, and a \leg\ unknot $K_0$ in $B$
with $tb(K_0)=tb(K')$, $rot(K_0)=rot(K')$.
By the classification of \leg\ unknots given by Eliashberg and Fraser
in \cite{eliashberg-fraser} $K_0$ is \leg\ isotopic to $K'$, thus also
to $K$. But $K_0$ is standard, while $K$ is non-loose: a contradiction.
\qed

\medskip
For a \lk\ $K$ the pair of integers $\bigl(rot(K),tb(K)\bigr)$
can be represented as a point in the plane (with integer coordinates 
of different parity).
Then the values of invariants allowed (by \ref{prop-antybenneq})
for non-loose {\bf un}knots are contained in one quadrant of 
the plane (dark shading in \ref{fig-nonloose};
the light shading denotes pairs $(t,r)$ satisfying
\ref{prop-bensw}, but not \ref{prop-antybenneq}). 
In fact, we only know the actual representants
(non-loose unknots with given rotation and Thurston-Bennequin 
invariant) for points $(rot, tb)$ 
on the boundary of this area (the black dots): 
those are the knots $\Gamma_{\pm1,\mp n}$  
(for any $n>0$) defined in \ref{ex-gammas}.

\midinsert
\centerline{\epsffile{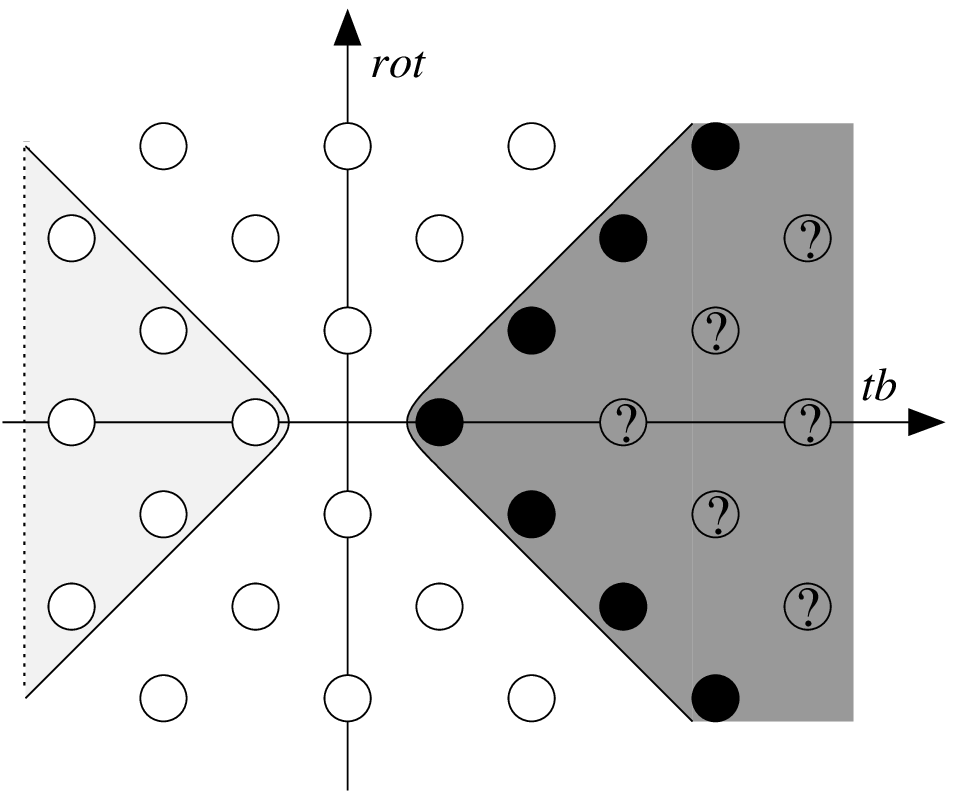}}
\figure{Non-loose unknots (known and potential).}
\definexref{fig-nonloose}{\the\figno}{figure}
\endinsert

\thm{Proposition}{The knots $\Gamma_{m,n}$ with $mn<0$ are non-loose.}
\definexref{prop-nonloosegammas}{\here}{propo}

\proof
The proof is analogous to the proof of \ref{propo-nonloose}. We show that
the complement of the knot $\Gamma_{m,n}$ is tight by embedding its
universal cover into the standard tight $\r^3$.

\smallskip
Let $\alpha=\arctan(m/n)$, so that $\pr_{\cal C}^{-1}(\Gamma_{m,n})$
lies in the plane $x=\alpha$. Then for $M=S^3\setminus \Gamma_{m,n}$ 
we have $M=M_f\cup M_b$, where $M_f=\pr_{\cal C}\bigl([0,\alpha]\times[0,1]^2\bigr)
\setminus\g mn$ and $M_b=\pr_{\cal C}\bigl([\alpha,1]\times[0,1]^2\bigr)
\setminus\g mn$. The universal cover $\widetilde{M_f}$ is diffeomorphic to 
$H_m\times\r$, where $H_m$ denotes a 2-dimensional disk with $m$ boundary 
points $p_1,\ldots,p_m$ removed (and similarly $\widetilde{M_b}\simeq H_n\times\r$;
the following argument applies to each of them).

The \cstr\ on $\widetilde{M_f}$ extends to a \cstr\ on
$\widetilde{M_f}\cup\{p_2\}\times\r\cup\ldots\cup\{p_m\}\times\r$,
isomorphic to the \cstr\ on $M_0$ used in the proof of
\ref{propo-nonloose}. 
Therefore, we have a contact embedding
$J:\widetilde{M_f}\to(\r^3,\zeta)$ with the boundary mapped
into the plane $x=0$ (actually onto, except for the $m-1$ lines
$\{x=0,z=z_i\}, i=2,\ldots,m\}$ which are ``images'' of
$\{p_i\}\times\r$). 

We want, however, to find an embedding whose image is contained in 
a proper subset of $\r^3$, so that there is enough space left
for other copies of $\widetilde{M_f}$ and $\widetilde{M_b}$.
It is well known that $(\r^3,\zeta)$ can be contactomorphically
embedded into an arbitrarily small open subset of itself.
Let $h: \r^3\to\{x^2+z^2< C\}$ be such a contact embedding.
We can assume that $h$ maps the plane $x=0$ into itself.

For our purposes, $h\circ J$ is better than $J$, but still not perfect,
and needs further fine-tuning. Let $R:\widetilde{M_f}\to\widetilde{M_f}$
be the rotation by the angle $2\pi/m$, cyclically permuting the $p_i$'s.
Then $h\circ J\circ(R\circ J^{-1}\circ h\circ J)^{m-1}$
embeds $\widetilde{M_f}$ into $H\times\r$, where $H$ is a subset
of the half-plane bounded by the line $x=0$ and $m$ arcs;
see \ref{fig-buildblock} with $m=3$.

\smallskip
\ref{fig-tree} shows the pattern of gluing together the pieces constructed 
in the described way. The flexibility in setting values of $J$ on
the boundary of $\widetilde{M_f}$ which we enjoy in the early stages
of this construction (cf.~the proof of Lemma 4.2.1 in \dy) is
sufficient to ensure that the embedding of $\widetilde{M_f}$ and
$\widetilde{M_b}$ glue together in a compatible way to give 
an embedding of the whole $\widetilde M$.

\bigskip
\rmk{Question} \definexref{que-exist}{\here}{question}
Let $r,t$ be integers satisfying the conditions $r+t\equiv 1 \pmod 2$,
$t\geq 3$, $|r|\geq t-3$. 
Does there exist a non-loose \leg\ unknot $K$ with $rot(K)=r$ and
$tb(K)=t$?
\bigskip

Of course, there do exist {\it loose} unknots for any given 
values of invariants (of different parity, of course).
This means that for some values there are (at least)
two non-isotopic knots. For example, we list all that 
is known about \leg\ unknots $K$ such that $rot(K)=0$
and $tb(K)=1$.

\midinsert
\centerline{\epsffile{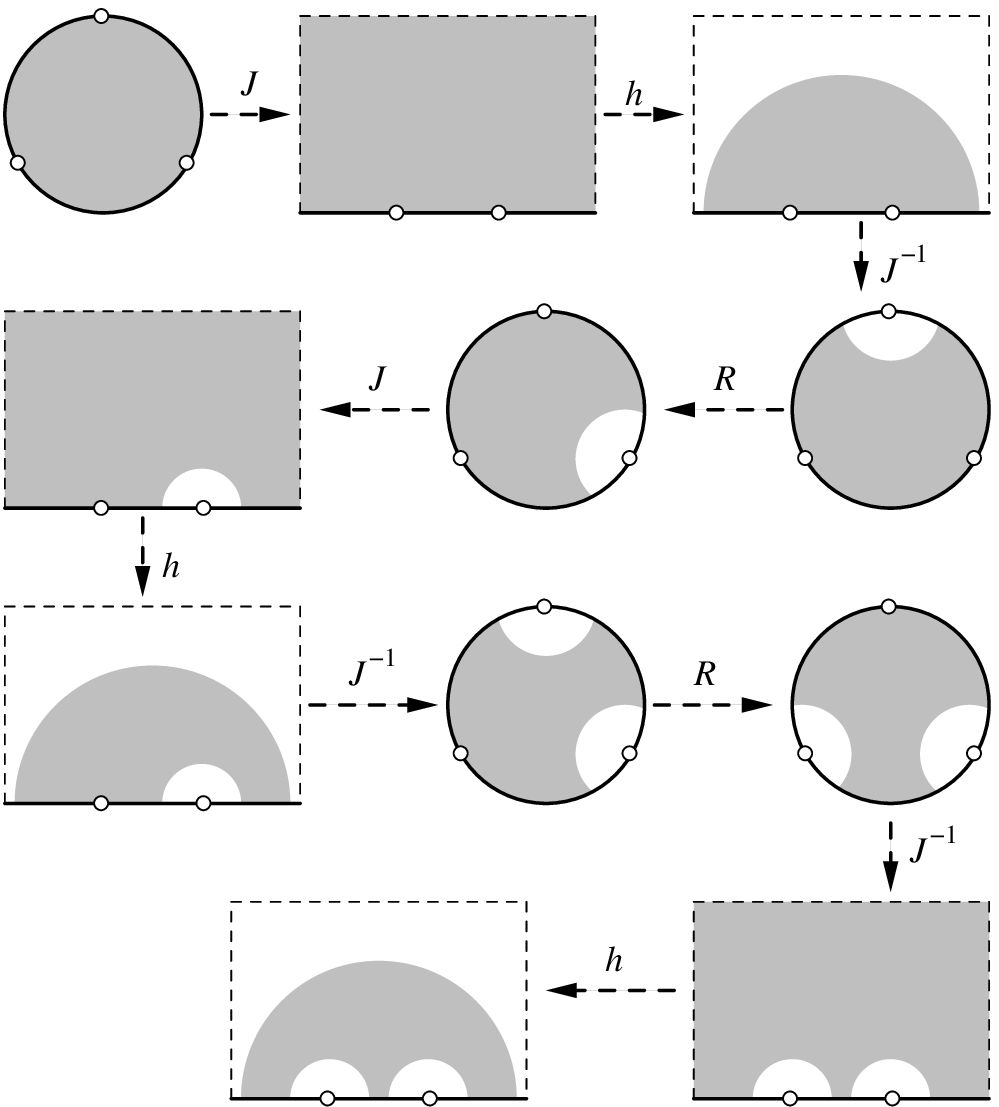}}
\figure{Constructing a building block for $\widetilde M \hookrightarrow \r^3$}
\definexref{fig-buildblock}{\the\figno}{figure}
\endinsert

\midinsert
\centerline{\epsffile{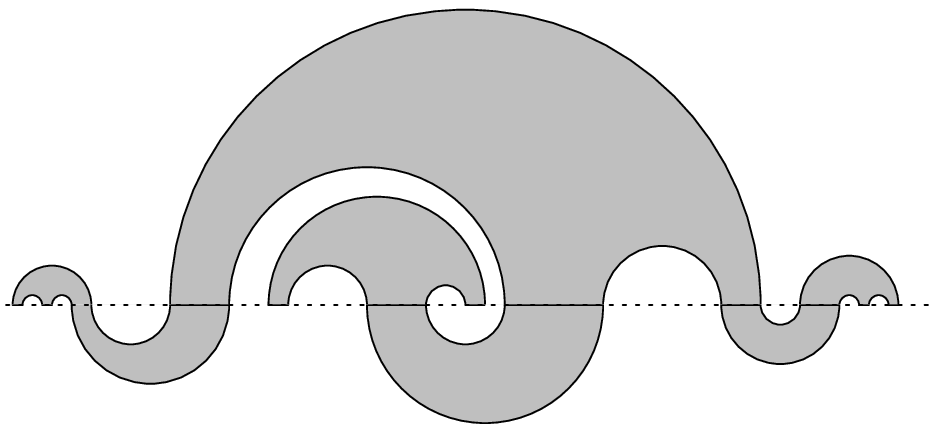}}
\figure{Gluing the building blocks together}
\definexref{fig-tree}{\the\figno}{figure}
\endinsert

\bigskip

\rmk{Example} \definexref{ex-kozuch}{\here}{exple}
Consider the four \leg\ unknots of \ref{fig-kozuch}
(each of them has $tb=1$ and $rot=0$):
\item{(i)} $\g 01\#\overline\g 01$
and $\g 10\#\overline\g 10$
are loose, while $\g 11$ and $\overline\g 11$
are not, thus neither $\g 11$ nor  $\overline\g 11$
is isotopic to $\g 01\#\overline\g 01$
or to $\g 10\#\overline\g 10$;
\item{(ii)} it is not known whether
$\g 01\#\overline\g 01$
and $\g 10\#\overline\g 10$
are isotopic;
\item{(iii)} it is not known whether
$\g 11$ and $\overline\g 11$
are isotopic;
\item{(iv)} it is not known if
there exist any other (non-isotopic) \leg\ unknots
with the same values of invariants.

\midinsert
\centerline{\epsffile{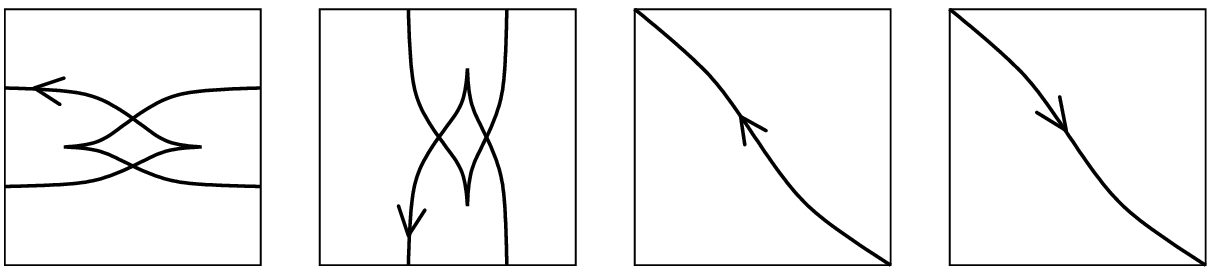}}
\hskip1cm $\g 01\#\overline\g 01$ \hfill $\g 10\#\overline\g 10$
\hfill $\g 11$ \hfill $\overline\g 11$ \hskip2cm \null
\figure{\lk s with $rot=0$ and $tb=1$.}
\definexref{fig-kozuch}{\the\figno}{figure}
\endinsert

\sec{References}

\bibliographystyle{siam}
\bibliography{contact,lanl}

\bye